\documentclass[a4paper,11pt, reqno]{amsart}

\usepackage[margin=2.75cm]{geometry}

\usepackage{amsmath}
\usepackage{booktabs}
\usepackage{amssymb}%$\nexists$
\usepackage{amsthm}%entorno para utilizar proof
\usepackage{color}
\usepackage{adjustbox}
\usepackage{graphicx}
\usepackage[ruled,vlined]{algorithm2e}
\usepackage{subcaption}%subfigure
\usepackage[authoryear]{natbib}
\usepackage{float}
\usepackage{hyperref}
\usepackage{multirow}
\usepackage{multicol}
\usepackage{soul}
\usepackage{algorithm2e}
\newtheorem{lem}{Lemma}
\newtheorem{dfn}{Definition}

\newtheorem{ex}{Example}

\def\A{\mathcal{A}}
\def\X{\mathcal{X}}

\newcommand{\dsum}{\displaystyle\sum}
\newcommand{\dmin}{\displaystyle\min}

\def\R{\mathbb{R}}
\def\B{\mathbb{B}}
\def\Z{\mathbb{Z}}
\def\cov{\mathcal{C}}

\usepackage{tikz}

\usepackage{pgfplots}
\pgfplotsset{compat=newest}
\usetikzlibrary{decorations.markings}

\definecolor{armygreen}{rgb}{0.19, 0.53, 0.43}

\definecolor{armygreen2}{rgb}{0.29, 0.13, 0.83}

\definecolor{Darkgreen}{rgb}{0,0.6,0}

\usepackage[colorinlistoftodos]{todonotes}
\definecolor{Granate}{rgb}{0.65, 0.04, 0.37}

\allowdisplaybreaks

\let\origmaketitle\maketitle
\def\maketitle{
	\begingroup
	\def\uppercasenonmath##1{} % this disables uppercasing title
	\let\MakeUppercase\relax % this disables uppercasing authors
	\origmaketitle
	\endgroup
}

\begin{document}
	
	\title[]{\Large Fairness in maximal covering location problems}
	
	\author[V. Blanco \MakeLowercase{and} R. G\'azquez]{{\large V\'ictor Blanco$^\dagger$ and Ricardo G\'azquez$^\dagger$}\medskip\\
		$^\dagger$Institute of Mathematics (IMAG), Universidad de Granada\\
		\texttt{vblanco@ugr.es}, \texttt{rgazquez@us.es}
	}
	
	\date{\today}
	
	\maketitle

		\begin{abstract}
This paper provides a general mathematical optimization based framework to incorporate fairness measures from the facilities' perspective to Discrete and Continuous Maximal Covering Location Problems. The main ingredients to construct a function measuring fairness in this problem are the use of: (1) ordered weighted averaging operators, a family of aggregation criteria very popular to solve multiobjective combinatorial optimization problems; and (2) $\alpha$-fairness operators which allow to generalize most of the equity measures. A general mathematical optimization model is derived which captures the notion of fairness in maximal covering location problems. The models are firstly formulated as mixed integer non-linear optimization problems for both the discrete and the continuous location spaces. Suitable mixed integer second order cone optimization reformulations are derived using geometric properties of the problem. Finally, the paper concludes with the results obtained on an extensive battery of computational experiments on real datasets. The obtained results support the convenience of the proposed approach.
		\end{abstract}
		

\section{Introduction}\label{sec:int}

The term \textit{fairness} is defined as ``the quality of treating people equally or in a way that is right or reasonable'' (Cambridge Dictionary). It is an abstract but widely studied concept in Decision Sciences in which some type of indivisible resources are to be shared among different agents. The importance of fairness issues in resource allocation problems has been recognized and well studied in a variety of settings with tons of applications in different fields \citep[see e.g.,][]{kelly1998rate,luss1999equitable,li2006fair,jiang2021rawlsian}. Fair allocations should imply impartiality, justice and equity in the allocation patterns, which are usually quantified by means of inequality measures that are minimized. Several measures have been proposed in the literature to this end, although the most popular one is the max-min (or min-max) approach which assures that the most \textit{harmed} agent in the share is as less damaged as possible \citep[see e.g.,][]{megiddo1974optimal,hayden1981voice,jaffe1981bottleneck,bertsekas1992data}. 
Other proposals of fairness measures are the minimum envy ~\citep{lipton2004approximately,caragiannis2009low,espejo2009comparison,netzer2016distributed} or certain families of ordered weighted averaging  criteria~\citep{ogryczak2006equity,hurkala2013fair,ogryczak2014fair}, among others.

In this paper, we analyze the notion of fairness in the context of Location Science. Facility location problems aim to \textit{optimally} determine the position of one or more facilities in order to satisfy the demand of a set of users. There is a vast variety of location problems, which are commonly classified by the nature of users and facilities, the optimization criteria, the type of demand to be satisfied, among many others \citep{hamacher1998classification}. The interested reader is referred to the recent book \citet{LocationScience2019} for a fresh view on the developments on both theoretical an applied aspects of location problems. 

We analyze here a particular family of location problems, the so-called Covering Location (CL) problems, that arise whenever a decision-maker aims to cover a given demand in case the facilities have a limited coverage area. These types of problems appear in many practical situations in which the services to be located are not able to satisfy the demand outside their coverage area. One of the most interesting CL problems is the Maximal Covering Location Problem (MCLP) whose goal is to find the positions of a given number of services, $p$, each of them with a coverage area, maximizing the covered demand of the users. This problem was firstly introduced by \citet{church1974maximal}, and since then the MCLP has attracted significant attention from both researchers and practitioners~\citep{murraywei}, both by its technical merit and practical interest. As in other Location Science problems, it is usual to analyze two different types of location spaces for the MCLP, in terms of the nature of the facilities that are to be located. In the \textit{discrete} MCLP, the services are selected from a given finite set of potential facilities, whereas in the \textit{continuous} MCLP, the services are allowed to be located at any position of the decision space (usually, the plane). 

This paper provides a general mathematical optimization based framework to incorporate fairness measures from the facilities' perspective to discrete and continuous MCLPs. We assume that a given number of services has to be located to maximize the demand covered by the different agents. In a fairly ideal solution, one would desire to ``independently'' maximize the covered demand of each of the services, not affecting negatively to the demand coverage of the others. However, since the demands are usually indivisible, in most cases, an advantageous solution for one service (one covering a large amount of the demand) harms others. As already occurs in other resource allocation problems, one may prefer to \textit{slightly} sacrifice the overall covered demand in order to balance the different covered demands among the open services. This might be applicable to the case of the location of public schools, in which it is preferable to find an homogeneous distribution of kids among the schools, or the location of routers with high capacities, where a ``good'' location for them would be the one in which the performance of all the routers can be better used instead of saturating some and leaving others covering a small amount of users. 

\subsection{Related work}

In many resource reallocation problems there has been special attention to the notion of fairness. For this reason, it has been widely recognized in the literature and, due also to the subjective implications the term fairness has, it has been studied in a wide variety of ways. The applications of fair resource allocation range from social or humanitarian contexts to engineering applications or location problems.

The social context is an important field where the fairness has been used. \citet{fairm2} highlights the importance of studying the notions of equity and justice and its incorporation into decision problems such as the fair allocation of water in irrigation farms. \citet{HUANG2019887} analyze the notion of fairness in the field of humanitarian logistics.

In engineering applications, as in telecommunication networks, where limited resources must be allocated among the different competing entities, the notion of fairness is also crucial~\citep[see e.g.,][]{Kleinberg1999,luss1999equitable,bonald2001,ogryczak2014fair}. 

The allocation of public resources is one of the most studied applications when analyzing equity and fairness. For instance, in the allocation of beds and ambulances in emergency medical services. In \citet[Chapter 4]{Leclerc2012}, the authors review  different measures previously proposed in the literature and study their incorporation to the allocation of ambulances. Another example of the importance of equity allocation comes from healthcare scheduling. Specifically, in the allocation of beds and other resources to patients where fairness allows to reduce mortality~\citep[see e.g,][]{8868095}.

In Location Science some equity measures have been also analyzed~\citep{espejo2009comparison,CHANTA2014228}. A review of measures for equitable facility location problems is given in \citet{MARSH19941,barbati2016equality}. 
	
Nevertheless, the notion of fairness can be interpreted in many different ways and there is no universal principle that is accepted as ``the most fair''~\citep{bertsimas2012efficiency}. The most usual strategy to incorporate fairness in resource allocation problems is to quantify the degree of fairness associated with any feasible action through a fairness measure, a function that maps each feasible solution into a real value indicating its degree of equity. Therefore, different measures have been proposed in the literature in order to capture the notion of fairness. The most popular one is the max-min (or min-max) approach which assures that the most \textit{damaged} agent in the share is as less damaged as possible \citep[see e.g.,][]{jiang2021rawlsian,kelly1998rate,li2006fair,luss1999equitable}. The max-min ratio is given by the maximum pairwise ratio between the resources allocated to the agents. Furthermore, some authors have incorporated fairness in decision problem through \emph{envy} measures, understood as the overall deviation in the gain (or cost) between all pairs of agents~\citep[see e.g.,][]{lipton2004approximately,caragiannis2009low,espejo2009comparison,CHANTA2014228,netzer2016distributed}. 

Other popular measure that has been widely investigated is the Jain's index~\citep{jainindex} which is obtained by normalizing the square mean of an entropy function. The Jain's index measures the \emph{equity} of an allocation pattern, being the Jain's index 1 for the most fair allocation in which all the agents receive (or pay) the same amount. This index was generalized by \citet{lan2010axiomatic}, where the authors propose a method to construct fairness schemes based on different axioms. The authors propose fairness measures generated by power functions generalizing the Jain's index.

An alternative approach that is very used in recent mathematical optimization based problems is to construct objective functions to guide the optimization approach to obtain fair solutions. This is the case of the so-called $\alpha$-fairness~\citep{kelly1998rate,mo2000fair,lan2010axiomatic,bertsimas2011price,bertsimas2012efficiency}. The $\alpha$-fairness functions were introduced by \citet{atkinson70} and was defined as the one maximizing the constant elasticity social welfare. This family of functions generalizes well-known criteria in the the allocation of resources, as the utilitarian principle ($\alpha=0$), which is neutral toward inequalities; proportional fairness ($\alpha=1$), that was introduced by \citet{nash1950}; and for $\alpha \rightarrow \infty$ the measure converges to the max-min criterion. As we will observe, other possibilities of $\alpha$ allow to find a trade-off between efficiency and fairness, opening room for different and suitable allocation patterns. 

An ideal approach to find solutions to fully satisfy all the agents involved in the allocation problem is to model it as a multiobjective problem. However, solving this problem may be computationally difficult in practice~\citep{ehrgott2000survey}. Moreover, the solution of a multiobjective problem is, in general, not a single solution but a solution set, the Pareto frontier, which can be useless in many decision problems. Thus, it is widely accepted to use aggregation operators of each of the agents utilities to find compromise solutions. Several aggregation criteria have been proposed to construct solutions satisfying all the different agents in the decision problem. The most popular aggregation criteria are: the weighted sum aggregation, the minimum of the  values (resulting in max-min problems) and the lexicographic order \citep{ogryczak2014fair}. 

Another general family of aggregation criteria that has been successfully applied is Ordered Weighted Averaging (OWA) functions, introduced by \citet{yager88}. These operators have been also applied to encourage fairness in the obtained solutions \citep[see e.g.,][]{ogryczak2006equity,hurkala2013fair,ogryczak2014fair}. An OWA operator maps a vector to a weighted sum of its coordinates in which the weights are assigned to their sorted values. These operators allow to generalize most of the robust statistical measures, as the mean, the minimum, the $k$-mean, the trimmed mean values or the Gini index, and have been applied in different fields~(see e.g., \cite{argyris2022fair,blanco2018locating,blanco2021multisource,marin2022soft}, among many others).

The use of OWA operators in Facility Location is not new and several authors have studied the incorporation of these operators to the classical problems through the so-called ordered median location problems \citep[see][]{puertofdez94}. Furthermore, when the weights defining the OWA function are nonnegative and monotone (non decreasing) the OWA operator is considered fair aggregation criteria, allowing to find fair distribution of the demand. 

In this paper, we analyze a novel version of one of the core family of problems in Facility Location, Covering Location problems. The location of the facilities in CL is characterized by the fact that the facilities are allowed to give service to the users at a limited  distance from them. This type of problems has been applied in different fields, as in the location of emergency services ~\citep{toregas1971location}, mail advertising~\citep{dwyer1981branch}, archaeology~\citep{bell1985location}, among many others. The interested reader is referred to \citet[Chapter 5]{LocationScience2019} for further details and recent advances on CL problems. In particular, we answer the question of how to incorporate fairness in the Maximal Covering Location problem (MCLP). In the MCLP, introduced by \citet{church1974maximal}, it is assumed the existence of a budget for opening facilities and the goal is to accommodate it to satisfy as much demand as possible. As usual in location problems, one can consider different frameworks based on the nature of the solution space for the facilities: discrete or continuous spaces. While the discrete setting is more adequate when locating  physical services, (as ATMs, stores, hospitals, etc), the continuous location space is known to be more adequate to determine the positions of routers, alarms or sensors, that can be more flexibly positioned. This type of space is also useful to determine the set of potential facilities that serves as input for a discrete version of the problem. One of the main difference between these two families of problems (from the mathematical optimization viewpoint) is that in the discrete one the distances between the facilities and the users are given as input data (or can be preprocessed before solving the problem), while in the continuous case, the distances are part of the decision and they must  incorporated to the optimization problem. 

As far as we know, the incorporation of fairness measures into CL problems has been studied in \citet{drezner2014maximin}, in which the max-min approach has been applied to the gradual maximal covering location problem. There, the authors incorporate the worst-case fairness criteria from the user's perspective, i.e. in order to enforce equity between the partial coverage of all users. In networks, \citet{rahmattalabi2020exploring} consider the selection of a subset of nodes to cover their adjacent nodes with fairness constraints with applications to social networks. \citet{Asudeh2020MaximizingCW} analyze a covering location problem with fairness constraints minimizing the pairwise deviations between the different covered sets. \citet{Korani2013TheHH} study a hub covering problem with \emph{equity} allocation constraints.

Despite of these applications, the efficiency measure used in the MCLP is the overall covered demand, that is, as much covered demand the better. However, when one looks at the individual utilities of each of the constructed facilities, one may obtain solutions with highly saturated facilities in contrast to others that only cover a small amount of demand, which results in unfair systems from the facilities' perspective. Moreover, in many situations this type of unfair solutions are also undesirable from the users' viewpoint which may see reduced the quality of the required service, as in the location of telecommunication servers which have a higher probability to fail in case of being saturated (being other capable to give service to these users) or in the student assignment process to schools, in which a higher number of alumni allocated to a school may deteriorate the education system. As far as we know, this problem has not been previously investigated in the literature in the context of Covering Location.

\subsection{Contributions}

In this paper, we provide a flexible mathematical optimization based framework to incorporate fairness measures from the facilities' perspective to Discrete and Continuous MCLPs. This generalization of the fairness measure for the MCLP is based on adequately combining the two main tools mentioned above: OWA operators and $\alpha$-fairness operators.

Our specific contributions in this paper are:
\begin{enumerate}
	\item To define a novel fairness measure combining OWA and $\alpha$-fairness operators that can be incorporated to the objective function of the MCLP.
	\item To describe a general mathematical optimization model which captures the notion of fairness from the facilities' perspective in MCLPs.
		\item To provide a mixed integer non-linear optimization formulation for the two main location spaces in the facility location problems: discrete and continuous. 
	\item To derive MISOCO reformulations for the problem, suitable to be solved with off-the-shelf optimization softwares.
	\item  To test on a battery of computational experiments the computational performance of the approaches and their managerial insights.
\end{enumerate}

\subsection{Paper structure}

The remainder of the paper is organized as follows. Section \ref{sec:preliminaries} is devoted to recall the used tools for our proposed method and states the notation for the rest of the paper. In Section \ref{sec:FOWA} we introduce the generalized fair maximal covering problem. In Section \ref{sec:formulations} we present a mathematical optimization formulation for the problem, both in the discrete and the continuous location spaces. The results of our computational experience are reported in Section \ref{sec:tests}. Finally, the paper ends with some conclusions and future research lines.

\section{Preliminaries}\label{sec:preliminaries}
In this section we introduce the notation used for the rest of the sections as well as the main results for analyzing the problem studied in this paper.

\subsection{Maximal Covering Location Problem}\label{sec:prel:MCLP}

Consider a finite set of demand points in a $d$-dimensional space, $\A = \{a_1, \ldots, a_n\} \subseteq \R^d$, indexed by the set $N=\{1, \ldots, n\}$. Each demand point $a_i \in \A$ has associated a non-negative demand weight $\omega_i$. Throughout the paper we often call a demand point interchangeably by the node $a_i$ or by the index $i$. Demand points may represent users or regions and the weights allow one to give more importance to different users or take into account the size/population of each of the regions. We are also given a potential set of facilities, $\X \subseteq \R^d$, where the \emph{services} are to be chosen. The set $\X$ is not necessarily finite. 

Each potential position for the facilities, $X \in \X$ is endowed with a coverage area. It is usual to define the coverage areas as Euclidean balls with certain \textit{coverage radii}. In this paper, we consider ball-shaped coverage areas in the form:
$$
\B_R(X)  =   \{z \in \R^d: \|z-X\| \leq R\}
$$
\noindent where $R>0$ is the radius and $\|\cdot\|$ is a $\ell_\tau$-based norm ($\tau \geq 1$), that is:
$$
\|x\|_\tau = \left(\dsum_{l=1}^d |x|^\tau\right)^{\frac{1}{\tau}},
$$
or a polyhedral norm with a symmetric polytope $B$ (with respect to the origin) in $\R^d$:
$$
\|x\|_B = \min\{\dsum_{g=1}^G |\beta_g|: x = \dsum_{g=1}^G \beta_g e_g\},
$$
where $\{\pm e_1, \ldots, \pm e_G\}$ are the extreme points of $B$.

A demand point $a\in \A$ is said to be \textit{covered by the facility} $X \in \mathcal{X}$ if it belongs to the coverage area of $X$, i.e., $a \in \B_R(X)$. 
Given $X \in \mathcal{X}$, we denote by $\cov(X)$, the index set of demand points covered by $X$, i.e., 
$$
\cov(X) = \{ i \in N: a_i \in \B_R(X)\}.
$$

The goal of the Maximal Covering Location Problem (MCLP), introduced by \citet{church1974maximal}, is to determine the positions of $p$ new facilities in $\X$, that is $\{X_1,\ldots,X_p\} \subseteq \X$, maximizing the ($\omega$-weighted) covered demand. Denoting  by $P = \{1, \ldots, p\}$ the index set for the facilities to be located, the MCLP can be formally formulated as the following optimization problem:
$$
\max_{X_1, \ldots, X_p \in \mathcal{X}} \dsum_{i \in \bigcup_{j\in P} \cov(X_j)} \omega_i.
$$
Two main families of MCLP arise based on the nature of the metric set from where the facilities are to be located. In case $\mathcal{X}$ is a finite set, one obtain the \textit{Discrete} MCLP that was introduced by \citet{church1974maximal}. On the other hand, in case $\mathcal{X}=\R^d$ one obtain the \textit{Continuous} MCLP that was firstly formulated by \citet{church1984} for the planar case ($d=2$). In this paper, we will analyze an unified and particular version of the MCLP for both type of location spaces.

%\subsection{Fairness}
\subsection{Ordered Weighting Averaging operators}\label{sec:prel:OWA}

As in most decision problems in which several agents are involved, the MCLP also exhibits a compromise in view of the users that require the services to be located or for the decision makers that construct and locate the facilities to give the different services to the users. In this paper, we analyze the MCLP from the point of view of the different decision makers that invest on the installation of the facilities to give service to the users. In an ideal solution one would desire to determine the position of the services to maximize, separately, the demand covered by each service. This approach will result on modeling the problem as a multiobjective mathematical optimization problem (with objectives the demands covered by each of the facilities) which can be cumbersome in practice. Instead, we aggregate the $p$ objective functions by means of Ordered Weighted Averaging (OWA) operators, introduced by \citet{yager88}. 

OWA operators are functions in the form $\Phi_\mathbf{\lambda}: \R^p \rightarrow \R$ with associated weighting vector $\mathbf{\lambda} = (\lambda_1, \ldots, \lambda_p)$, with $\lambda_j \in [0,1],\; \forall j \in \{1,\ldots,p\}$ and $\sum_{j = 1}^{p} \lambda_j = 1$. We denote as $P = \{1, \ldots, p\}$ the index set for coordinates of vectors in $\R^p$. For $\mathbf{W} = (W_1, \ldots, W_p) \in \R^p$, the OWA operator is defined as:
\begin{equation} \label{eq:owafunction}\tag{\rm OWA}
	\Phi_\mathbf{\lambda}(W_1,\ldots, W_p) = \dsum_{j = 1}^{p} \lambda_j W_{(j)},
\end{equation}
\noindent where $W_{(j)}$ is the $j$th-largest value in the vector $\mathbf{W}$, i.e., $W_{(j)} \in \{W_1, \ldots, W_p\}$ such that $W_{(1)} \leq \ldots \leq W_{(p)}$. OWA operators are weighted sums of the different criteria, but where weights are associated with the position of the criteria when they are sorted in non decreasing order.

When aggregating different criteria by means of an OWA operator, a crucial step is to provide the adequate weights $\lambda$ that are assigned to the sorted sequence.  The most popular OWA operators ($\min$ -- $\lambda=(1,0, \ldots,0)$ and $\max$-- $\lambda=(0,\ldots,0,1)$) serve as reference weights to define the notion of \textit{orness} of a vector of $\lambda$-weights defining an OWA operator: 
$$
{\rm orness}(\lambda) = \dsum_{j=1}^p \frac{p-i}{p-1} \lambda_i
$$

The degree of orness emphasizes the higher (better) values or the lower (worse) values in a set of attributes associated with the different agents/services. Given a vector of $\lambda$ weights, as closer its orness to $1$, closer to the $\min$-operator while as closer to $0$, closer to the $\max$-operator. Assuming that all the criteria are to be minimized, the $\min$-operator allows one to generate solutions protected under worst-case scenario (pessimistic), while the $\max$-operator produces solution in which the best situation for all the criteria is assumed (optimistic). In between, one can find an equilibrium between those extreme choices. In particular, for $\lambda=(\frac{1}{p}, \ldots, \frac{1}{p})$-- the mean operator, its orness degree takes value $0.5$. In Table \ref{table:owas} we show a list with some of the most popular OWA operators and their orness degree.

\begin{table}[H]
\centering
\adjustbox{width=\textwidth}{\begin{tabular}{llll}\hline
$OWA$ & $\lambda$ & Operator & ${\rm orness}$\\\hline
Average & $\lambda_ j = \frac{1}{p}$ & $\frac{1}{p} \dsum_{j\in P} W_j$ & $\frac{1}{2}$\\
Minimum & $\lambda_1=1, \lambda_j=0$ ($j\geq 2$) & $\min_{j\in P} W_j$ & 1\\
$k$-Average & $\lambda_j=\frac{1}{k}$ ($j\leq k$), $\lambda_j=0$ ($j>k$) & $\frac{1}{k} \dsum_{j=1}^k  W_{(k)}$ & $1 - \frac{k-1}{2(p-1)}$\\
$\alpha$-Min-Average & $\lambda_1=\frac{1}{1+(p-1)\alpha}$, $\lambda_j=\frac{\alpha}{1+(p-1)\alpha}$ ($j\geq 2$) & $\frac{1-\alpha}{1+(p-1)\alpha} \min_{j\in P} W_j + \frac{\alpha}{1+(p-1)\alpha}  \dsum_{j\in P}  W_{j}$ & $\frac{-p\alpha+p+2\alpha}{2p\alpha-2\alpha+2}$\\
Gini & $\lambda_j = \frac{2(p-j)+1}{p^2}$  for all $j$& $\frac{1}{p^2} \dsum_{j\in P} W_j +\frac{2}{p^2} \dsum_{j\in P} (p-j) W_{(j)}$ & $\frac{4p+1}{6p}$\\
Harmonic & $\lambda_j =\frac{1}{p} \left(H(p)-H(j-1)\right)$  ($H(k) = \dsum_{\ell = 1}^{k} \frac{1}{\ell}$) & $\frac{1}{p} \dsum_{j\in P} \left(H(p)-H(j-1)\right)  W_{(j)}$ & $\frac{3}{4}$\\
\hline
\end{tabular}}
\caption{Some examples of fair OWA operators and its orness.\label{table:owas}}
\end{table} 

It is clear that an OWA operator (identified with a $\lambda$-weight) is not uniquely determined by its orness degree (unless its orness degree is in $\{0,1\}$ or $p=2$). Thus, several optimization-based methods have been proposed in order to construct, with different paradigms, $\lambda$-weights with a given orness degree $\beta \in (0,1)$ \citep[see, e.g.,][]{filev1995analytic,fuller2001analytic,liu2004properties}. The main idea when searching for $\lambda$-weights with a given orness degree is to solve problems in the form:
\begin{align*}
\min &\;\; \mathcal{L}(\lambda)\\
\mbox{s.t. } &{\rm orness}(\lambda) = \beta, \\
& \lambda \in \R^p_+,
%& \lambda_1 \geq \ldots \geq \lambda_p \geq 0
\end{align*}
where $\mathcal{L}$ is a loss function measuring some properties of the weights. For instance, if $\mathcal{L}(\lambda) =-\sum_{j=1}^p \lambda_j \log \lambda_j$ one obtain the maximal entropy monotone OWA \citep{o1988aggregating}, or choosing $\mathcal{L}(\lambda) = \sum_{j=1}^p (\lambda_j - \overline{\lambda})^2$ one obtains the minimum variance weights \citep{fuller2003obtaining}, where $\overline{\lambda}$ stands for the mean of the vector $\lambda$.

There are several ways to sort a list of (unknown) values and compute the OWA operator. The most known are the OT proposed by \citet{ogryczak03} and the BEP reformulation proposed by \citet{BEP14}. In the former, the authors provide a suitable linear optimization representation of the problem of minimizing the sum of the $k$ largest (equivalently, smallest) linear functions on a polyhedral set in $\R^d$. This representation is extended to the minimization of monotone OWA functions by means of a telescopic sum of $k$-sum functions. 

\subsection{Fairness} \label{sec:fairness}

As already mentioned, fairness is a crucial requisite when making decisions in many situations in which different agents are involved. Since the notion of fairness is by itself weak and it is not clear, in general, a list of axioms/desirable properties have been stated for a vector of values to be considered fair. \citet{lan2010axiomatic} provides an axiomatic approach to fairness measures which served at starting point to derive new fairness measures fulfilling a list of minimum required properties. Here, we provide some useful of desirable properties of fairness schemes. Let be $\Psi: \R^p_+ \rightarrow \R$ a fairness scheme, and $\mathbf{W} = (W_1,\ldots,W_p) \in \R^p_+$ the allocation of resources. Some of these properties were stated in \citep{lan2010axiomatic,jainindex,barbati2016equality}.
\begin{enumerate}
	\item {\bf Continuity:} $\Psi$ is a continuous function. This axiom assures that, locally, small changes in the allocation do not significantly affect the measure.
	\item {\bf Population size independence}: Equal resource allocations are, eventually, independent of the number of users, i.e., $\lim_{p\to \infty} \frac{\Psi(\mathbf{1}_{p+1})}{\Psi(\mathbf{1}_{p})} =1$.
	\item {\bf Pareto optimality:} If $W_j \leq \bar{W}_j, \forall j \in \{1,\ldots, p\}$, and $W_j < \bar{W}_j$ for at least some $j$, then $\Psi(\mathbf{W}) \leq \Psi(\mathbf{\bar{W}})$.
\item {\bf Symmetry:} $\Psi(W_1,\ldots,W_p) = \Psi(W_{\sigma(1)},\ldots, W_{\sigma(p)}),
	$ where $\sigma$ is an arbitrary permutation of the indices.	
	\item {\bf Bounded:} The value of allocation given by the scheme is bounded.
	\item {\bf Scale and metric independence:} The measure is independent of scale, i.e, the unit of measurement does not matter.
	
\end{enumerate} 

Here, we consider a particular measure of fairness which has gained attention in the operations research community~\citep{bertsimas2011price,bertsimas2012efficiency}, the $\alpha$-fairness scheme, that was introduced by \citet{atkinson70}. According to this scheme, the decision maker decides on the allocation of the resources by maximizing the constant elasticity social welfare operator $\Psi_\alpha$ parameterized by $\alpha \geq 0$ and for $\mathbf{W} = (W_1,\ldots,W_p) \in \R^p_{+}$, the amount of resources allocated into $p$ agents,
\begin{equation*}%\label{eq:alphafairness}
	\Psi_\alpha(W_1, \ldots, W_p) = \left\{\begin{array}{cl}
	\frac{1}{1-\alpha} \dsum_{j =1}^p W_j^{1-\alpha} & \mbox{if $\alpha \geq 0$, $\alpha\neq 1$,}\\
	\dsum_{j=1}^p \log(W_j) & \mbox{if $\alpha=1$}.
	\end{array}\right.
\end{equation*}

The parameter $\alpha$ is known as the \emph{inequality aversion parameter} since it controls the difference between the allocation resources $W_j$. Consider an agent $k$ with lower amount of resources $W_k$ than the agent $k'$ with amount of resources $W_{k'}$. If we increase the resources in $W_k$ then we would have a higher welfare than if we increase in $W_{k'}$. Thus, an increase in the resources of $W_k$ would be more desirable to reduce the unfairness in the allocation of the resources. This can be performed by tuning, adequately, the parameter $\alpha$ since when it increases, the difference between the values of resources for the agents decreases \citep[for a constructive proof of this, see ][]{lan2010axiomatic}. This property, also known Principle of transfer or Pigou-Dalton \citep{Erkut1993}, yielding then to fair solutions.

Therefore, moving the $\alpha$ value, one can represent classical fairness measures.  For $\alpha \rightarrow \infty$, the scheme converges to the max-min fairness, whereas if $\alpha=0$ one obtain the sum, which is neutral toward inequalities. For $\alpha=1$, one obtain the proportional fairness measure that has been widely studied in the literature \citep{nash1950}.

\section{The generalized Fair Maximal Covering Location Problem} \label{sec:FOWA}

As already mentioned, the Maximal Covering Location Problem (MCLP), in its different versions, can be seen as a resource allocation problem, in which the overall demand of the \textit{covered} users is shared among the different services that are located. Thus, a high coverage of the total demand (no matter which service is providing the coverage) is appropriate from a global perspective, but from an individual viewpoint, one can easily get unfair allocation patterns. Furthermore, the MCLP usually exhibits multiple optimal solutions, that is, different subsets of $p$ services covering the maximum possible covered demand, and then, optimization solvers output an arbitrary one, possibly not the fairest. 

In the following example we illustrate this situation in a toy instance.
\begin{ex}
	We consider a randomly generated set of $200$ demand points with coordinates in $[0,20]\times [0,20]$. We assign to each of them a random integer weight in $[0,100]$. We assume that three $(p=3)$ services are to be located chosen from the same set of demand points. The coverage areas for all the potential location of the services are disks of radius $5$. The solution obtained by Gurobi for the classical MCLP is shown in Figure \ref{fig:ex0a}. In such a solution, $65.7\%$ of the demand is covered, and the distribution of (weighted) users among the services is $(1723, 2365, 2804)$, that is, there are two services covering close to one thousand more clients than the other one. Other feasible (non optimal for the MCLP) solution for the problem is also shown in Figure \ref{fig:ex0b}, in which $62.58\%$ of the demand is covered and whose distribution of covered demand is $(2126, 2162, 2278)$. This solution, although covers $3\%$ less demand than the classical MCLP, is clearly much more balanced than the MCLP solution, since all the users cover approximately the same demand, but still efficient.
	
	\begin{figure}[H]
		\begin{subfigure}[b]{.48\linewidth}
			\centering \fbox{{\centering \includegraphics{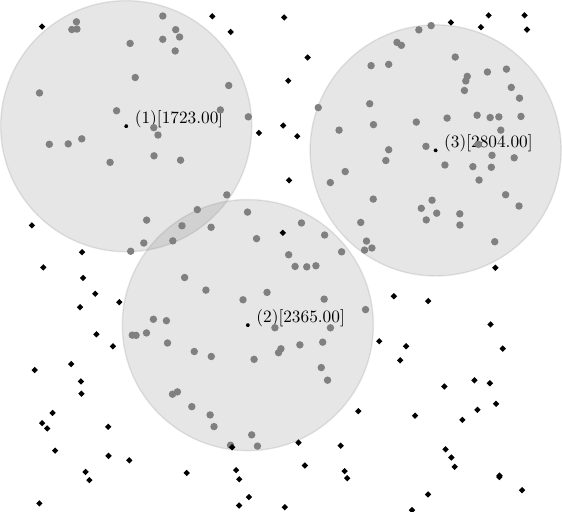}}}
			\caption{Classical MCLP.}\label{fig:ex0a}
		\end{subfigure}\hfill
		\begin{subfigure}[b]{.48\linewidth}
			\centering \fbox{{\centering \includegraphics{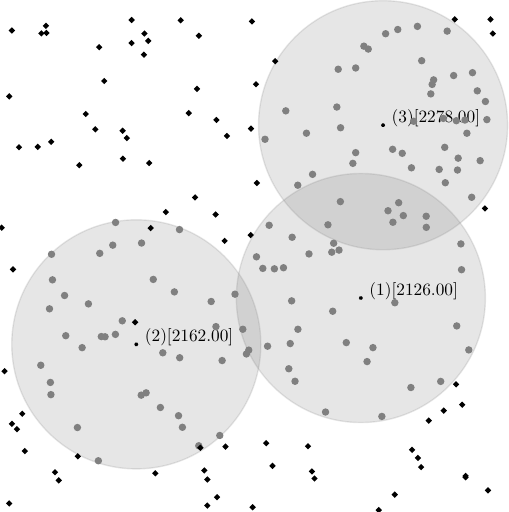}}}
			\caption{Fairer solution for MCLP.}\label{fig:ex0b}
		\end{subfigure}
		%\fbox{\input{fig2Mean.tex}}~\fbox{\input{fig2Min.tex}}
		\caption{Optimal solutions obtained with the MCLP.\label{fig:ex0}}
	\end{figure}
\end{ex}

In what follows we introduce a new fairness measure that we incorporate into the MCLP in order to provide fair coverage of the demands.
\begin{dfn}
	Let $\alpha \geq 0$ and $\lambda \in \R_+^p$ with $\sum_{j=1}^p \lambda_j=1$ and $\lambda_1 \geq \cdots \geq \lambda_p$. The $(\alpha,\mathbf{\lambda})$-fair operator is a function $F_{\alpha,\lambda}: \R^p_+ \rightarrow \R^+$ defined as:
	\begin{equation}\label{eq:fowa}\tag{$(\alpha,\mathbf{\lambda})$-Fair}
	F_{\alpha, \mathbf{\lambda}} (W_1, \ldots, W_p) : = \left\{\begin{array}{cl}
	\frac{1}{1-\alpha} \dsum_{j =1}^p \lambda_j W_{(j)}^{1-\alpha} & \mbox{if $\alpha\neq 1$,}\\
	\dsum_{j=1}^p \lambda_j \log(W_{(j)}) & \mbox{if $\alpha=1$}.
	\end{array}\right.
	\end{equation}
\end{dfn}

It is straightforward to check that $F_{\alpha, \mathbf{\lambda}}$ is a concave function that verifies all of the desirable properties described in Section \ref{sec:fairness}. 

The \eqref{eq:fowa} operator depends of $p+1$ parameters ($\alpha$ and the $\lambda$-weights) and combines the $\alpha$-fairness measure introduced by \cite{atkinson70} and the OWA operators introduced in \cite{yager88}. In case $\alpha=0$, the operator turns into the $\mathbf{\lambda}$-OWA operator, whereas if  $\mathbf{\lambda}=(\frac{1}{p}, \ldots, \frac{1}{p})$ it becomes the $\frac{1}{p}$-weighted  $\alpha$-fairness function $\frac{1}{p}\; \Psi_\alpha$. This combination of these two operators allows us to derive a unified framework to deal with fairness in maximal covering location problems, that we detail below.

We are given a set of demand points in a $d$-dimensional space, $\A = \{a_1, \ldots, a_n\} \subseteq \R^d$, indexed by the set $N=\{1, \ldots, n\}$, and each of the points endowed with a demand weight $\omega_i \geq 0$ for all $i\in N$. We are also given a metric space $\mathcal{X} \subseteq \R^d$ endowed with a distance measure $\|\cdot\|$ and a radius, $R$, for each $X\in \mathcal{X}$, which is assumed to be the same for all the facilities to be located (although it is not a limitation of the results provided in this paper). The goal of the Fair MCLP is to find the position of the $p$ facilities to locate, $X_1, \ldots, X_p \in \mathcal{X}$, maximizing the $(\alpha,\lambda)$-fair operator of the demands covered by these services. Formally, given $\mathbf{\lambda}$ and $\alpha$, the $(\alpha,\mathbf{\lambda})$-Fair Maximal Covering Location Problem ($(\alpha,\lambda)$-FMCLP, for short) can be model as the following optimization problem: 
\begin{equation}\label{fmclp:0}\tag{${\rm FMCLP}_{\alpha,\mathbf{\lambda}}$}
F^*_{\alpha, \mathbf{\lambda}} = \max_{X_1, \ldots, X_p \in \X} F_{\alpha,\mathbf{\lambda}} (W(X_1), \ldots, W(X_p))
\end{equation}
\noindent where $W(X_j)$ is the covered demand of facility $X_j$, assumed that each demand point is accounted as covered by at most one of the facilities. We denote by $\mathcal{W}_{\alpha;\mathbf{\lambda}} =$ $(W(X_1^*)$, $\ldots,$ $W(X_p^*)) \in \R^p_+$ a coverage vector of each of the facilities in the problem above. We also denote by $\mathcal{W}^{\rm sum}_{\alpha;\mathbf{\lambda}} = \sum_{j=1}^p W(X_p^*)$ the total covered demand in the solution and by $\mathcal{W}^{\rm min}_{\alpha;\mathbf{\lambda}} = \dmin_{j=1, \ldots, p} W(X_p^*)$ the demand covered by the service with smallest coverage in the solution.

The problem that we introduced above allows one to determine the position of the facilities that are fair from the individual viewpoint, but, how much is one willing to lose when imposing fairness to the MCLP? \emph{The price of fairness} in any allocation rule is a notion studied in \cite{bertsimas2011price} in order to measure the efficiency loss under a fair allocation compared to the one that maximizes the overall sum of the users utilities. In our case, the solution of \eqref{fmclp:0} is compared against  the solution of the classical MCLP in order to know how far is a fair solution to the best coverage of the give demand.
\begin{dfn}
	The price of fairness measure for \eqref{fmclp:0} is defined as the index:
	\begin{equation}\label{eq:pof}\tag{PoF}
	{\rm PoF} \eqref{fmclp:0} = \dfrac{\mathcal{W}^{\rm sum}_{0;(\frac{1}{p}, \ldots, \frac{1}{p})}-\mathcal{W}^{\rm sum}_{\alpha;\mathbf{\lambda}}}{\mathcal{W}^{\rm sum}_{0;(\frac{1}{p}, \ldots, \frac{1}{p})}},
	\end{equation}
\end{dfn}

The price of fairness indicates the relative deviation of the covered demand when solving the maximal $(\alpha,\lambda)$-FMCLP with respect to the solution of the classical MCLP which attains the maximal possible coverage. Thus, the price of fairness is a value between $0$ and $1$ measuring how close is the effectiveness of the obtained fair solution with respect to the most effective covering solution. A price of fairness equal to $0$ indicates that \eqref{fmclp:0} is able to construct a fair allocation without loss of efficiency (at the maximum possible coverage). In contrast, a price of fairness with value $1$ indicates that \eqref{fmclp:0} got the worst global coverage. Thus, as closer this measure to $0$ the better. In general it provides the percent loss of coverage with respect to the maximal possible coverage of an instance,  allowing one to quantify the price one has to pay when imposing $(\alpha,\lambda)$-fairness. 

On the other hand, one can also measure how far is a fair solution from the fairest share, which is obtained when solving the max-min covering location problem, i.e., comparing the demand covered (in our fair MCLP) by the service covering the smallest demand in the solution with respect to the solution in which the coverage of the service covering the least demand is maximized. This measure was called in  \cite{bertsimas2012efficiency} the \emph{price of efficiency}.

\begin{dfn}
	The price of efficiency measure for \eqref{fmclp:0} is defined as:
	\begin{equation}\label{eq:poe}\tag{PoE}
	{\rm PoE}  \eqref{fmclp:0} = \dfrac{\mathcal{W}^{\rm min}_{0;(1,0, \ldots, 0)}-\mathcal{W}^{\rm min}_{\alpha;\mathbf{\lambda}}}{\mathcal{W}^{\rm min}_{0;(1,0, \ldots, 0)}},
	\end{equation}
\end{dfn}

The Price of Efficiency is interpreted as the percent loss in the minimum demand coverage guarantee compared to the maximum minimum covered demand guarantee.  This index also takes value in $[0,1]$, in which a value of $0$ means that the $(\alpha,\lambda)$-FMCLP obtains the fairest solution, while a value of $1$ indicates the least fair solution in which there is a service not covering any demand.

Other widely used measure of fairness is the \emph{envy}. In the MCLP the envy of facility positioned in $X_j\in \mathcal{X}$, whose covered demand is $W(X_j)$, for a facility located at $X_k\in \mathcal{X}$, whose covered demand is $W(X_k)$, is defined as:
$$
{\rm envy}(X_j,X_k) = \max\{0, W(X_k)-W(X_j)\}
$$
that is, facility $j$ suffers envy of value $W(X_k)-W(X_j)$ from facility $k$ in case $k$ covers more demand than $j$. The overall envy of a set of $p$ facilities $X_1, \ldots, X_p\in \mathcal{X}$ is the overall sum of all the pairwise envies. From this, the \emph{Gini index} is defined as the ratio of this total envy and the all the covered demand by the facilities multiplied by $2p$ (the overall number of pairwise comparisons):
\begin{dfn}
	The Gini index is defined as:
\begin{equation}\label{eq:gini}\tag{Gini}
{\rm Gini}(X_1, \ldots, X_p) = \frac{ \dsum_{j, k \in P} {\rm envy}(X_j,X_k) }{2p \dsum_{j\in P} W(X_j)}
\end{equation}
\end{dfn}

We will see in our computational experience that the family of $(\alpha,\mathbf{\lambda})$-FMCLP exhibits differences when varying the values of $\alpha$ and $\lambda$ with respect to the three measures that we described  above (PoF, PoE and Gini). A trade-off solution between these measures will be desirable from the point of view of efficiency and also of fairness.

\section{Mathematical optimization formulations for $(\alpha,\mathbf{\lambda})$-FMCLP}\label{sec:formulations}

In this section we derive a suitable mathematical optimization formulation to model the $(\alpha,\mathbf{\lambda})$-FMCLP. We will present different formulations for the problem for both the discrete case ($\X$ being a finite pre-specified set) and the continuous case ($\X = \R^d$). The nature of the domain of this problem directly affect the development of resolution strategies for it. 

A general formulation for the problem considers the following decision variables:
$$
x_{ik} = \begin{cases}
1, & \text{if node $a_i$ is covered by the $k$-th selected facility in $\X$},\\
0, & \text{otherwise},\end{cases} \quad \text{ for all } i \in N,\: k \in P.
$$
and $X_k \in \R^d$: coordinates of the $k$-th selected facility in $\X$.

\eqref{fmclp:0} can be formulated as follows:
\begin{subequations}\label{fmclp:1}
	\begin{align}
	\max\; & F_{\alpha,\mathbf{\lambda}}\Big( \dsum_{i \in N} \omega_i x_{i1}, \ldots,  \dsum_{i \in N} \omega_i x_{ip}\Big) \\
	\mbox{s.t. } &  \dsum_{k \in P} x_{ik} \leq 1,\; \forall i \in N, \label{eq:fmclp2}\\
	& a_i \in \B(X_k)  \mbox{ if $x_{ik}=1$},\; \forall i \in N, k \in P, \label{eq:fmclp3}\\
	& x_{ik} \in \{0,1\},\; \forall i \in N, k \in P, \\
	& X_k \in \X,\; \forall k \in P, 
	\end{align}
\end{subequations}
The objective function aims to maximize the  $(\alpha,\mathbf{\lambda})$-fairness of the demand coverage by each of the facilities. Constraint \eqref{eq:fmclp2} assures that each covered demand point is counted at most once as covered. \eqref{eq:fmclp3} ensures the adequate definition of the $x$-variables. 

The feasible set of the problem above will described then by a set of linear and second-order cone representable inequalities on binary and continuous variables. One of the main difficulties of the model above stems on the representation of the objective function $F_{\alpha, \lambda}$ which consists of the following two ingredients:
\begin{description}
	\item[Sorting:] Representing the order given by the OWA operator into an optimization problem is a difficult challenge. In Section \ref{sec:prel:OWA} we describe the two most popular formulations for this operator on the values $W_k^{1-\alpha}$ for $k\in P$. A third representation is based on the $x$-variables that we consider in our problem, by sorting the selected facilities in $\X$ in non-increasing order of the demand coverage, i.e, enforcing the following constraints:
	\begin{equation}\label{sorting}\tag{Sorting}
	\dsum_{i \in N} \omega_i x_{ik} \leq \dsum_{i \in N} \omega_i x_{ik+1},\; \forall k \in P,
	\end{equation}
	\item[$(1-\alpha)$-powers:] Observe that the $(1-\alpha)$ powers of the coverage of each facility appear in the objective function. Denoting by $W_k = \dsum_{i \in N} \omega_i x_{ik}$ and by $Z_k = W_{k}^{1-\alpha}$, for $k \in P$, assuming that $W_1 \leq \ldots \leq W_p$, the objective function above can be written as the linear function:
	$$
	F_{\alpha,\mathbf{\lambda}}(W_1, \ldots, W_p) = \frac{1}{1-\alpha} \dsum_{k\in P} \lambda_k Z_k
	$$
	as long as it is fulfilled that $Z_k \leq W_k^{1-\alpha}$ (for $\alpha<1$) or $Z_k \geq W_k^{1-\alpha}$ (for $\alpha>1$) for all $k\in P$. Assuming that $\alpha$ is rational, we get that there exists $p, q \in \Z_+$ with $p \geq q \geq 1$ and $\gcd(p,q)=1$ such that:
	$$
	\frac{1}{1-\alpha} = \begin{cases} 
	\frac{p}{q} & \mbox{if $\alpha<1$},\\
	-\frac{q}{p} & \mbox{if $\alpha>1$}
	\end{cases}
	$$
	Thus, the power-constraints can be rewritten as:
	\begin{equation}\label{powers}\tag{Powers}
	Z_k^p \leq W_k^q, \text{for all $k\in P$}.
	\end{equation}
	
	These constraints can be conveniently rewritten as a set of quadratic second-order cone constraints following a simplification of the results in \cite{BEP14}. % In Algorithm \ref{alg:SOCctrs} we detail how to efficiently rewrite these constraints in a suitable way.
\end{description}

In the rest of the section we describe how to represent constraints \eqref{eq:fmclp3} in a suitable mathematical optimization formulation. This representation highly depends on the nature of the set of potential coordinates for the facilities $\X$. We analyze the cases in which $\X$ is a finite set and the one where $\X=\R^d$.

\subsection{Continuous facilities} \label{sec:FOWA:c}

We analyze here the case where the potential set from which the coordinates of the services are chosen is the entire space, i.e., $\X = \R^d$.  In this case, the norm-based \textit{covering contraints} \eqref{eq:fmclp3} can be rewritten as
\begin{equation}\label{eq:norms}\tag{Norms}
\|X_k - a_i\| \leq R + U_i (1-x_{ik}),  \forall i \in N, k \in P.
\end{equation}
where $U_i$ is a big enough constant ($U_i>\|a_i-a_{i'}\|$ for all $i'\in N$). It ensures that in case $i$ is allocated to the $k$th selected facility ($X_k$), then $a_i$ must belong to $\B_{R}(X_k)$. 

In case $\|\cdot\|$ is the Euclidean norm, it is well known that these constraints can be re-written as the following set of linear and second-order cone inequalities:
\begin{subequations}
\begin{align}
v_{ikl} \geq a_{il} - X_{kl}, & \: \forall i \in N, \: \forall k \in P,\label{norm:1}\\
v_{ikl} \geq -a_{il} + X_{kl},& \: \forall i \in N, \: \forall k \in P, \label{norm:2}\\
s_{ik} \geq \dsum_{l=1}^d v_{ikl}^2, & \: \forall i \in N, \: \forall k \in P,\label{norm:3}\\
s_{ik} \leq R + U_i\,(1-x_{ik}), & \: \forall i \in N, \: \forall k \in P, \label{norm:4}\\
v_{ikl} \geq 0, & \: \forall i \in N, \: \forall k \in P,\\
s_{ik} \geq 0, & \: \forall i \in N, \: \forall k \in P.
\end{align}
\end{subequations}
where $v_{ikl}$ represents (via \eqref{norm:1} and \eqref{norm:2}) the absolute value $|a_{il}-X_{kl}|$ (here $a_{il}$ and $X^t_{kl}$ stand for the $l$-th coordinate of demand point $a_i$ and the $k$-th facility, respectively) and $s_{ik}$ determines the euclidean distance (by \eqref{norm:3}) as the squared sum of the absolute differences between the coordinates of $a_i$ and $X_{k}$. Consequently, the problem simplifies to a  mixed-integer second-order cone optimization (MISOCO) problem, which can be solved using  any of the available off-the-shelf solver.

Actually, by the results in \cite{BEP14}, one can also consider block-norm based distances (deriving linear optimization models) or $\ell_\tau$-norms (with $\tau\geq 1$) inducing again mixed-integer second-order cone optimization problems. 

As already observed in other maximal covering location problems (see \cite{BG21,BGS21}) the norm-based constraints \eqref{eq:norms} can also be rewritten as linear constraints. This \textit{linearization} is based on projecting out the $X$-variables by ensuring that these can be constructed from the $x$-variables. The observation that will allow this formulation is a direct consequence of the following lemma whose proof is straightforward:
\begin{lem}\label{lem:1}
	Let $N_1,\dots,N_{p} \subseteq N$ be $p$ nonempty disjoint subsets of $N$. 
	Then, there exists $X_1, \ldots, X_p \in \R^d$ such that the points of $N_k$ are covered by $X_k$, for all $k \in P$,  if and only if
	$$
	\bigcap_{i \in N_k} \B_{R}(a_i) \neq \emptyset, \forall k \in P.
	$$
\end{lem}
From this result, \eqref{eq:norms} can be replaced by:
\begin{equation}\label{eq:fmclp:linear1} 
\dsum_{i \in Q} x_{ik} \leq |Q|-1, \forall k \in P,\: \forall Q \subseteq N: \bigcap_{i \in Q} \B_{R}(a_i) = \emptyset, \\
\end{equation}
These exponentially many constraints assure that those subsets of \textit{incompatible} demand points cannot be covered by the same facility. Thus, if a solution in the $x$-variables, $\bar x$, fullfilling \eqref{eq:fmclp:linear1}, we can construct the coordinates of the services that must verify
$$
X_k  \in \bigcap_{i \in N : \atop \bar x_{ik}=1} \B_{R}(a_i), \text{for all $k \in P$. }
$$
Actually, these coordinates can be found  (in poly-time) for a given feasible solution $\bar x$ either solving a second order cone optimization problem or by solving a one-center facility location problem. This strategy has been successfully applied in the literature by \citet{BG21,BGS21}.

Nevertheless, the above formulation requires incorporating the exponentially many constraints \eqref{eq:fmclp:linear1}. However, they can be simplified reducing from exponential to polynomial many constraints by means of Helly's Theorem \citep{Helly:1923} \citep[see also][]{DanzerGruenbaumKlee:1963}. Invoking that result, only intersections of $(d+1)$-wise balls are needed to check, allowing to replace Constraints \eqref{eq:fmclp:linear1} by:
\begin{equation}\label{eq:fmclp:linear} 
\dsum_{i \in Q} x_{ik} \leq |Q|-1, \forall k \in P,\: \forall Q \subseteq N: \bigcap_{i \in Q} \B_{R}(a_i) = \emptyset, \text{ with } |Q|=d+1. \\
\end{equation}

Despite this simplification, the number of constraints is still $O(n^d)$ making the problem  difficult to solve. Then, we propose to consider an incomplete formulation (removing \eqref{eq:fmclp:linear}) and iteratively incorporating these constraints \textit{on-the-fly}, as needed.

The selection of the constraints to incorporate in each iteration is found using the following separation strategy:  After solving the problem with none or part of the constraints \eqref{eq:fmclp:linear} a solution, say $\bar{\mathbf{x}}$, is obtained.
Then, for each $k \in P$ the define set $Q_k= \{i \in N: \bar x_{ik} =1\}$. 
One can check for the validity of the set $Q_k$ as a feasible cluster of demand points for our problem by solving the $1$-center problem for the points in such a set.
In case the optimal coverage radius obtained is less than or equal to $R$, one can conclude that $Q_k$ is a valid subset of demand points that can be covered by the same server. 
Otherwise, the solution violates the relaxed constraints, and thus we add the cut
\begin{equation}\label{Qk}
\dsum_{i \in Q_k} x_{ik'} \leq |Q_k|-1, \forall k' \in P,
\end{equation}
to ensure that such a solution is no further deemed feasible and thus obtained again. 

The $1$-center problem with Euclidean distances on the plane is known to be solvable in polynomial time \citep[see e.g.,][]{EH72}. Extensions to higher dimensional spaces and generalized covering shapes have been recently proved to be also poly-time solvable \citep{BP21}.

\subsection{Discrete facilities} \label{sec:FOWA:d}

Let us assume that the potential set of facilities if finite, that is, $\X = \{b_1,\ldots,b_m\} \subseteq \R^d$. We denote by $M = \{1,\ldots,m\}$ its index set. The model in this case can be simplified taking into account that the subset of potential facilities that are able to cover each single demand point can be pre-computed. It allows also to avoid the use of the $X$-variables, replacing them by the following decision variables to decide which of the potential facilities from $\{b_1,\ldots,b_m\}$ are open and which is the position of the demand covered by each facility in the ordered vector.
$$
y_{jk} = \begin{cases} 1, & \text{if the covered demand of facility $j$ is the $k$-th largest},\\
0, & \text{otherwise},\end{cases} \quad \text{ for all } j \in M,\: k \in P.
$$
Then, \eqref{fmclp:1} can be alternatively formulated as:
\begin{subequations}
%	\makeatletter
%	\def\@currentlabel{${\rm FMCLP}^{\rm D}$}
%	\makeatother
%	\label{fmclp:D}
%	\renewcommand{\theequation}{${{\rm FMCLP}_{\arabic{equation}}^{\rm D}}$}
	\begin{align}
	\max\; & \frac{1}{1-\alpha}\dsum_{k \in P} \lambda_k Z_k \\
	\mbox{s.t. } & \eqref{eq:fmclp2}, \eqref{sorting}, \eqref{powers} \label{eq:generals} \\
	& x_{ik} \leq \dsum_{j \in \mathcal{C}(i)} y_{jk},\; \forall i \in N, k \in P, \label{eq:fmcl:d:1} \\
	& \dsum_{j \in M} y_{jk} = 1,\; \forall k \in P, \label{eq:fmclp:d:2}\\
	& \dsum_{k \in P} y_{jk} \leq 1,\; \forall j \in M, \label{eq:fmclp:d:3}\\
	& x_{ik} \in \{0,1\},\; \forall i \in N, \forall k \in P, \\
	& y_{jk} \in \{0,1\},\; \forall j \in M, k \in P. 
	\end{align}
\end{subequations} 

Apart from \eqref{eq:generals}, the \textit{covering} constraints  \eqref{eq:fmclp3} are rewritten using the $y$-variables using  \eqref{eq:fmcl:d:1}--\eqref{eq:fmclp:d:3}.  Constraints \eqref{eq:fmcl:d:1} assure that the demand points can be assigned to a facility if it is sorted in any position (equivalently, if it is open). Constraints \eqref{eq:fmclp:d:2} enforces that a single facility is assigned to a position and \eqref{eq:fmclp:d:3} that at most one position is assigned to a facility (those facilities not assigned to a positions will be not open). Both constraints together with \eqref{sorting} assure that exactly $p$ facilities are open, each of them in a different order in the coverage sequence. 

The classical formulations for the MCLP use one-index binary variables to determine the \emph{open} facilities. Note that in our case it is not enough, since the positions of the activated services in the sorting coverage sequence are needed to allow the allocation of the demand points to a facility only in case it is open.

\section{Computational study} \label{sec:tests}

In this section, we report the results of the computational experience that we have run in order to analyze the performance of the proposed methoddologies.

We consider the dataset provided in \cite{IRSloc} where the coordinates (latitude and longitude) of Residential Schools and student hostels operated by the federal government in Canada are detailed. The coordinates of this dataset have been normalized in the unit square. 

From the whole dataset we have constructed different instances with sizes, $n$, ranging in $\{45,90,120,179\}$ (the first $n$ demand nodes in the complete instance).  The demands have been randomly generated by a uniform distribution in $(0,1)$. The set of generated instances are available in the github repository \href{https://github.com/vblancoOR/fmclp}{github.com/vblancoOR/fmclp}.
The number of facilities to be located, $p$, ranges in $\{5,10,15,20\}$, and we consider the radius for the services, $R$, ranging in $\{0.1,0.15\}$. 

We run our models by choosing the $\lambda$-parameters of the OWA operator in  $\{\texttt{W},\texttt{C},\texttt{K},\texttt{D},\texttt{G},\texttt{H}\}$ (see Table \ref{table:owas}, where \texttt{C} - Minimum, \texttt{D} - $\alpha-$Min-Average, \texttt{G} - Gini, \texttt{H} - Harmonic, \texttt{K} - $k-$Average,  and \texttt{W} - Average). The $\alpha$-parameter was chosen in $\{0,0.5,1,2\}$ except for problems under the $\texttt{C}$-objective, where the value of $\alpha$ does not affect the result, and we only solved the  instances for $\alpha=0$. 

Finally, we have solved all these instances for the two types of domains: discrete and continuous (we refer them in the rest of the currently section as \texttt{Disc} and \texttt{Cont}). With that, a total of 672 instances for each type of location space (a total of 1344 instances) have been solved in our computational study.
For the continuous MCLP, we observed that the nonlinear formulation \eqref{eq:norms} have a worse performance than the linear formulations proposed in the literature \citep[see][]{BG21,BGS21}, as expected. Thus, we provide the results of using formulation \eqref{Qk} in our experiments.

All the experiments have been run on a virtual machine in a physical server equipped with 8 threads from a processor AMD EPYC 7402P 24-Core Processor, 64 Gb of RAM and running a 64-bit Linux operating system. 
The models were coded in Python 3.7 and we used Gurobi 9.1 as optimization solver. A time limit of 2 hour was fixed for all the instances.

We analyze in this section the computational performance of our approach as well as the managerial insights of the obtained solutions in terms of fairness and efficiency.

\subsection{Computational performance}

In order to analyze the computational difficulty of the proposed approach, in Table \ref{t:times} we summarize the obtained results in terms of consumed CPU time, MIP Gap and unsolved intances. In column \texttt{Time} we show the average required computational times for those instances solved up to optimality within the time limit (2 hours). In column \texttt{GAP0} we report the percent of instances that were solved in the time limit. The rest of the columns show the information for those instances not optimally solved. \texttt{GAP1} shows the percentage of unsolved instances but MIPGap less than $1\%$, whereas  \texttt{GAP5} is the average percentage of instances that are with MIPGap in $(1\%,5\%]$. Finally in column \texttt{GAP unsolved}, we show the average MIPGap of the instances that we are not able to solve within the time limit. All the results are averaged for each fixed value of $n$ (number of demand points), $p$ (number of facilities to locate) and the type of location space (\texttt{Disc} and \texttt{Cont}).

The first observation is that the continuous problems are more computationally challenging to solve than the discrete ones. Even being both problems non linear in most of the cases because of the $\alpha$-fairness objective function, the distance-based constraints needed for the continuous problems make those problems much more difficult to solve (in terms of CPU time and unsolved instances).

The table also highlights the percentage of solved instances. As expected, as the number of facilities increases, the number of solved instances decreases. The MIPGap, however, are relatively small, being around 75\% the instances that were solved within the time limit, with less than 5\% of MIPGaps, 55\% being with less than 1\% MIPGap. The unique exception are the instances with $n=45$ and $p=20$. In those instances, a $100\%$-coverage can be obtained in all cases, but the number of feasible allocation is huge, being the underlying combinatorial problem very hard to solve.

\begin{table}[h]
	\centering
	\adjustbox{max width=\textwidth}{\begin{tabular}{llrrrrrrrrrr}
		\hline
		$n$ & $p$ & \multicolumn{2}{c}{\texttt{CPUTime (secs.)}} & \multicolumn{2}{c}{\texttt{GAP0}} & \multicolumn{2}{c}{\texttt{GAP1}} & \multicolumn{2}{c}{\texttt{GAP5}} & \multicolumn{2}{c}{\texttt{GAP unsolved}}\\
		\cmidrule(lr){3-4}\cmidrule(lr){5-6}\cmidrule(lr){7-8}
		\cmidrule(lr){9-10} \cmidrule(lr){11-12}
		    &     & \multicolumn{1}{c}{\texttt{Disc}} & \multicolumn{1}{c}{\texttt{Cont}} & \multicolumn{1}{c}{\texttt{Disc}} & \multicolumn{1}{c}{\texttt{Cont}} & \multicolumn{1}{c}{\texttt{Disc}} & \multicolumn{1}{c}{\texttt{Cont}} & \multicolumn{1}{c}{\texttt{Disc}} & \multicolumn{1}{c}{\texttt{Cont}} & \multicolumn{1}{c}{\texttt{Disc}} & \multicolumn{1}{c}{\texttt{Cont}} \\
		\hline		
		\multirow{4}{*}{45} & 5 & 2.00 & 25.50 & 100\% & 100\% & 100\% & 100\% & 100\% & 100\% & 0\% & 0\% \\
		 & 10 & 1672.71 & 1557.52 & 71.43\% & 42.86\% & 80.95\% & 71.43\% & 85.71\% & 83.33\% & 9.44\% & 4.70\% \\
		 & 15 & 192.04 & 159.60 & 11.90\% & 11.90\% & 28.57\% & 38.10\% & 54.76\% & 78.57\% & 44.68\% & 6.09\% \\
		 & 20 & 66.41 & 524.13 & 9.52\% & 9.52\% & 19.05\% & 23.81\% & 30.95\% & 57.14\% & 67.15\% & 60.52\% \\
		\hline
		\multirow{4}{*}{90} & 5 & 142.73 & 329.32 & 100\% & 80.95\% & 100\% & 100\% & 100\% & 100\% & 0\% & 0.03\% \\
		 & 10 & 366.22 & 2003.43 & 61.90\% & 33.33\% & 97.62\% & 83.33\% & 100\% & 97.62\% & 0.24\% & 0.98\% \\
		 & 15 & 1174.07 & 191.58 & 21.43\% & 11.90\% & 69.05\% & 35.71\% & 80.95\% & 92.86\% & 3.13\% & 2.10\% \\
		 & 20 & 56.66 & 1487.77 & 9.52\% & 7.14\% & 21.43\% & 16.67\% & 40.48\% & 38.10\% & 14.85\% & 7.74\% \\
		\hline
		\multirow{4}{*}{120} & 5 & 41.27 & 850.91 & 85.71\% & 83.33\% & 100\% & 100\% & 100\% & 100\% & 0.09\% & 0.17\% \\
		 & 10 & 1269.24 & 275.26 & 59.52\% & 9.52\% & 95.24\% & 50.00\% & 97.62\% & 88.1\% & 0.52\% & 2.36\% \\
		 & 15 & 958.53 & 868.72 & 30.95\% & 14.29\% & 78.57\% & 21.43\% & 85.71\% & 80.95\% & 2.05\% & 3.44\% \\
		 & 20 & 51.28 & 445.71 & 9.52\% & 9.52\% & 30.95\% & 11.90\% & 66.67\% & 64.29\% & 7.96\% & 4.47\% \\
		\hline
		\multirow{4}{*}{179} & 5 & 102.51 & 2306.03 & 78.57\% & 57.14\% & 100\% & 69.05\% & 100\% & 100\% & 0.03\% & 1.73\% \\
		 & 10 & 1897.44 & \texttt{TL} & 38.10\% & 0\% & 78.57\% & 9.52\% & 100\% & 28.57\% & 0.79\% & 10.85\% \\
		 & 15 & 1642.58 & 82.43 & 7.14\% & 2.38\% & 28.57\% & 9.52\% & 73.81\% & 45.24\% & 3.26\% & 23.16\% \\
		 & 20 & 2224.03 & 17.19 & 7.14\% & 4.76\% & 14.29\% & 7.14\% & 45.24\% & 42.86\% & 6.45\% & 7.38\% \\		
		\hline				
	\end{tabular}}
	\caption{Average CPU time (\texttt{CPUTime}), percentage of instances solved within time limit (\texttt{GAP0}), percentage of instances solved with less than 1\% and 5\% of gap (\texttt{GAP0} and \texttt{GAP5}), and average of MIP GAP of unsolved instances (\texttt{GAP unsolved}) averaged by size $n$ and number of centers $p$ for the discrete location problem.}
	\label{t:times}
\end{table}

\begin{figure}[H]
	\begin{center}
			\includegraphics{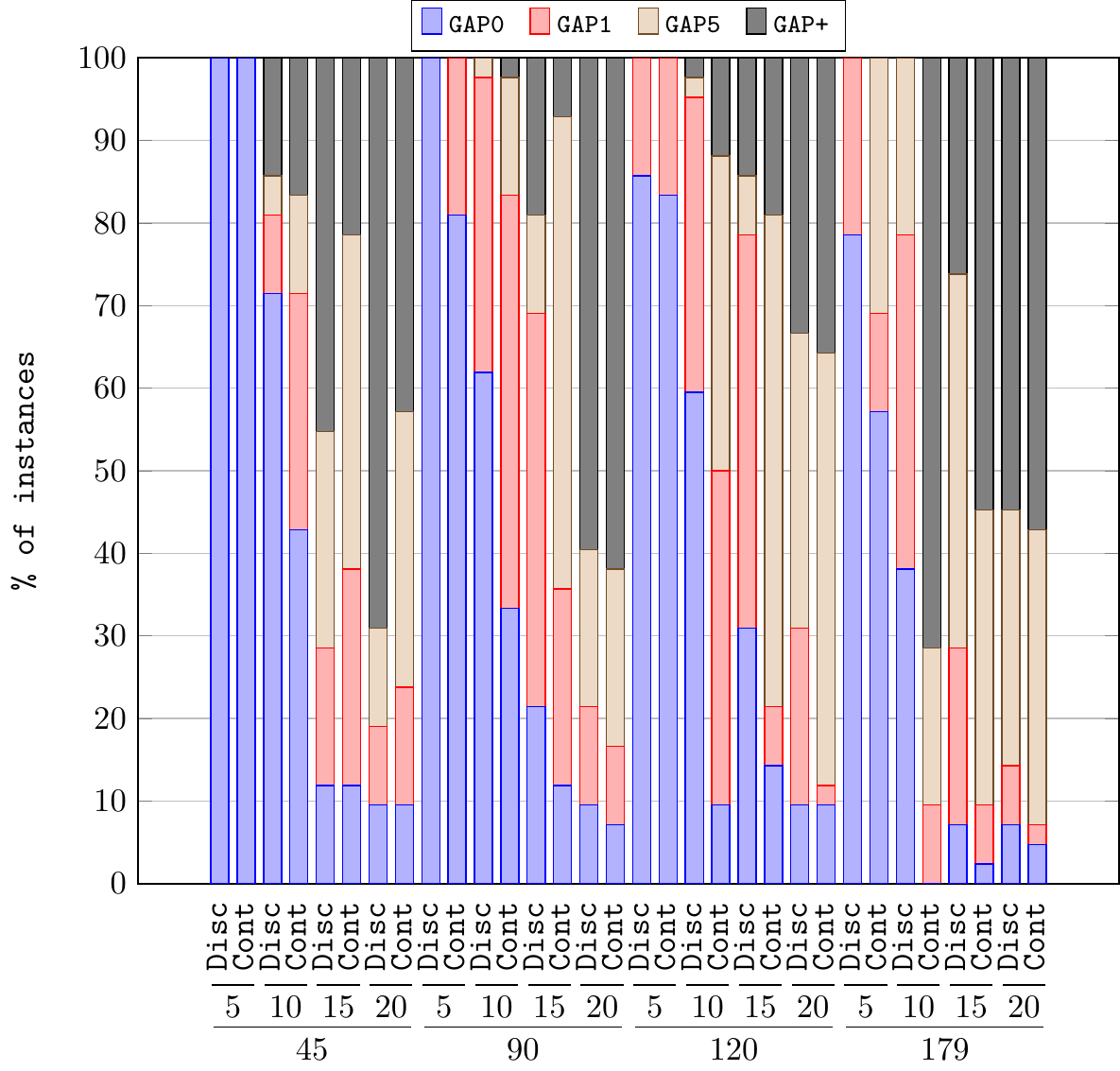}	
	\end{center}
\caption{Percentage of instances optimally solved within the time limit (\texttt{GAP0}), unsolved but with MIPGap less than $1\%$ (\texttt{GAP1}), unsolved by with MIPGap in $(1\%,5\%]$ (\texttt{GAP5}) and unsolved with MIPGap greater than $5\%$ (\texttt{GAP+}) averaged for each fixed value of $n$ (number of demand points), $p$ (number of facilities to locate) and the type of location space (\texttt{Disc} and \texttt{Cont}).}
	\label{f:stackedbar-table}
\end{figure}

In Figure \ref{f:stackedbar-table} we show the results of columns \texttt{GAP0}, \texttt{GAP1} and \texttt{GAP5}. Each bar represents all the instances for each combination of $n$, $p$ and Disc/Cont location problem. Each bar is partitioned into the four possible sets: optimally solved within the time limit (\texttt{GAP0}), unsolved but with MIPGap less than $1\%$ (\texttt{GAP1}), unsolved by with MIPGap in $(1\%,5\%]$ (\texttt{GAP5}) and unsolved with MIPGap greater than $5\%$ (\texttt{GAP+}).

Although percentage of solved instances within the time limit decreases as $p$ increases, the percentage of instances that have less than a 5\% gap is stable for each value of $n$. In addition we observed for values of $p\leq 10$, we were able to solve in most of the cases with less than one 5\% for the two types of location spaces.

\begin{figure}[H]
	\begin{center}
	\begin{subfigure}[b]{.5\linewidth}
		\includegraphics[scale=1]{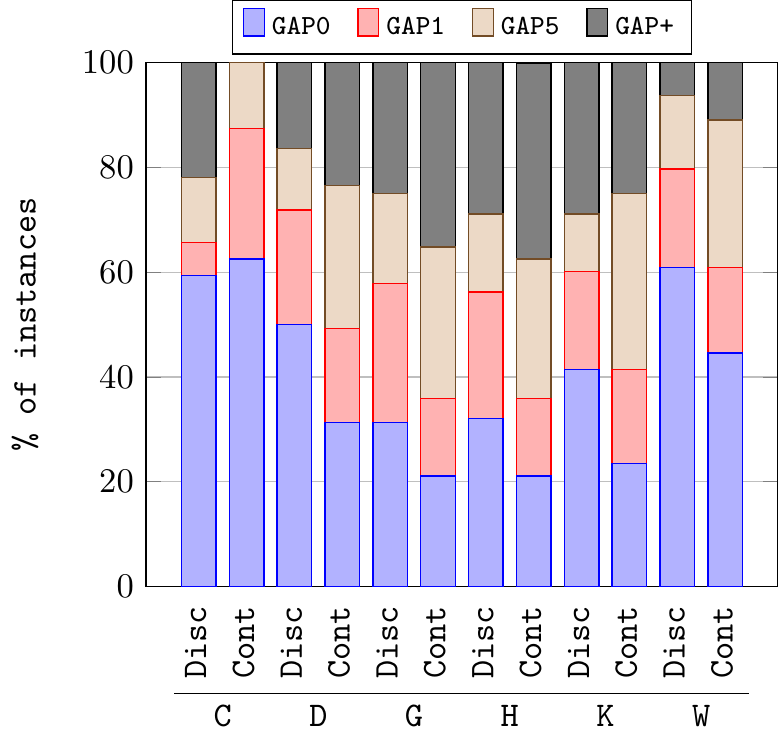}	
		\caption{By type of problem $\lambda$.}
	\end{subfigure}~
\begin{subfigure}[b]{.5\linewidth}
	\includegraphics[scale=1]{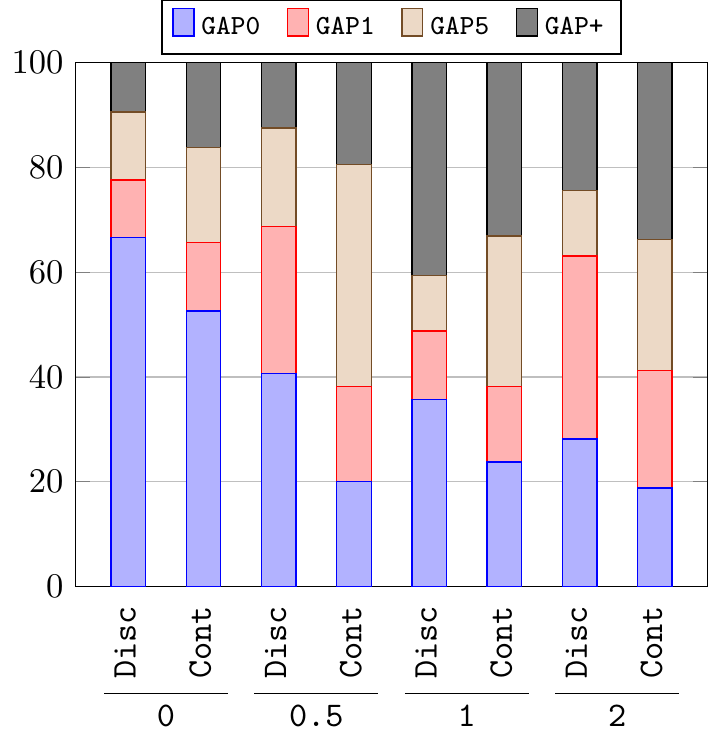}
	\caption{By value of $\alpha$.}
\end{subfigure}
\end{center}
	\caption{Percentage of instances optimally solved within the time limit (\texttt{GAP0}), unsolved but with MIPGap less than $1\%$ (\texttt{GAP1}), unsolved by with MIPGap in $(1\%,5\%]$ (\texttt{GAP5}) and unsolved with MIPGap greater than $5\%$ (\texttt{GAP+}) averaged for parameters involved in the problem $\lambda$--problem type and $\alpha$.}
	\label{f:stackedbar}
\end{figure}

In Figure \ref{f:stackedbar} we show a similar table but averaging by the type of problem ($\lambda$) and by the values of $\alpha$.
We observe that \texttt{C} and \texttt{W} were less computationally demanding than the others, with close to $60\%$ of the instances optimally solved. We also observed that around $70\%$ of the instances are solved with less than a 5\% gap.

However, for the value of $\alpha = 1$ the figure shows a higher percentage of instances with gap greater than 5\%, caused by the piecewise linear representation of Gurobi of the logarithmic function, which requires a new set of variables and constraints in the problem.

\subsection{Managerial insights}

In what follows we analyze the quality of the obtained solutions  with the different approaches in terms of fairness and efficiency.  As noted above, the MIPGaps for the instances unsolved within the time limit is small, being most of the instances assumed to be optimally solved (with certain degree of accuracy). In order to show the quality of the solutions, we consider those instances with a MIPGap smaller than 5\%. 

\begin{figure}[H]
	\begin{subfigure}[b]{.45\linewidth}
		\centering 
		\includegraphics[scale=0.9]{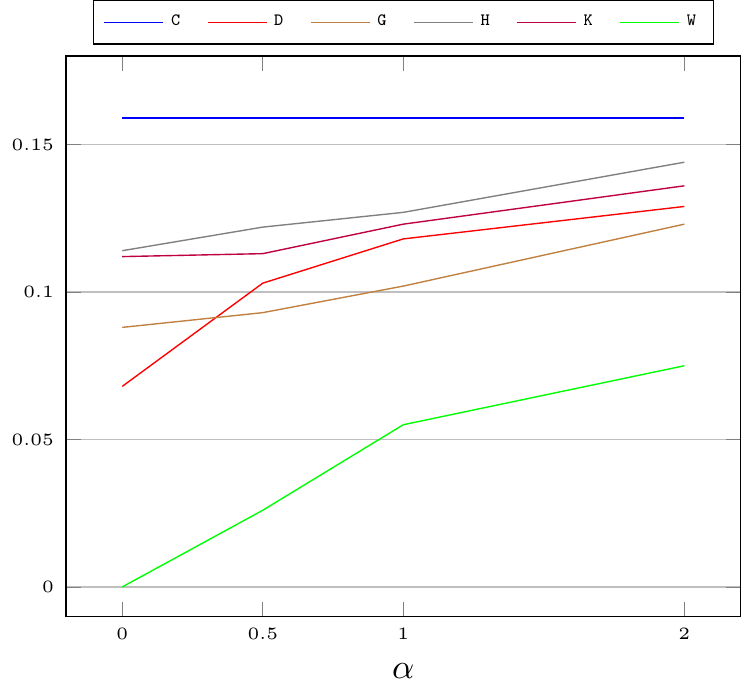}
		\caption{Price of fairness (PoF).}
	\end{subfigure}~\begin{subfigure}[b]{.45\linewidth}
		\centering 
		\includegraphics[scale=0.9]{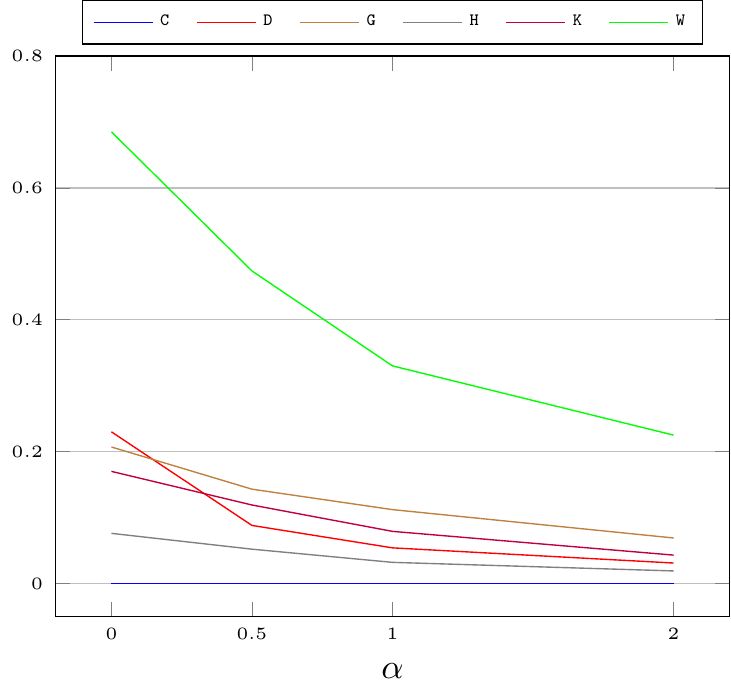}
		\caption{Price of efficiency (PoE).}
	\end{subfigure}
	\caption{PoF and PoE averaged for each value of $\alpha$ for \texttt{Discrete} location problem.}
	\label{f:discrete_np}
\end{figure}

\begin{figure}[H]
	\begin{subfigure}[b]{.45\linewidth}
		\centering 
		\includegraphics[scale=0.95]{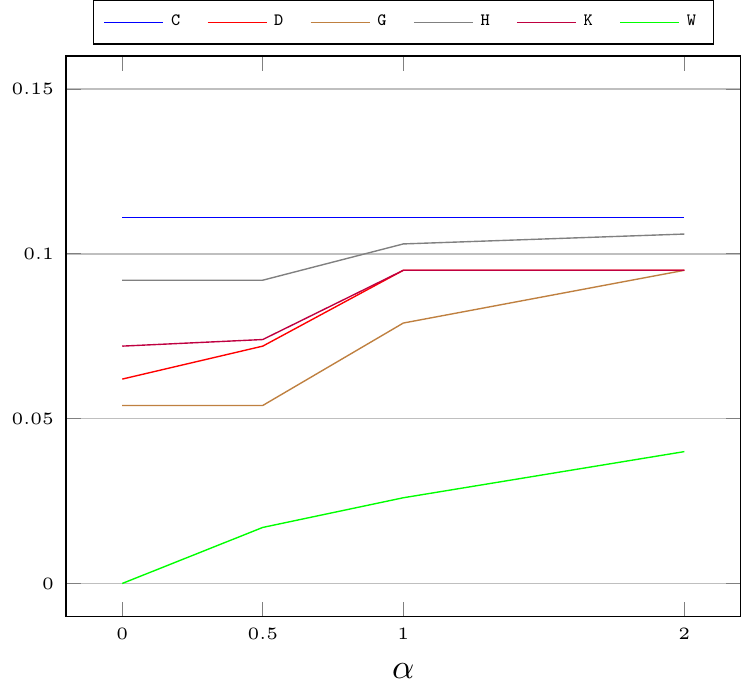}
		\caption{Price of fairness (PoF).}
	\end{subfigure}~\begin{subfigure}[b]{.45\linewidth}
		\centering 
		\includegraphics[scale=0.95]{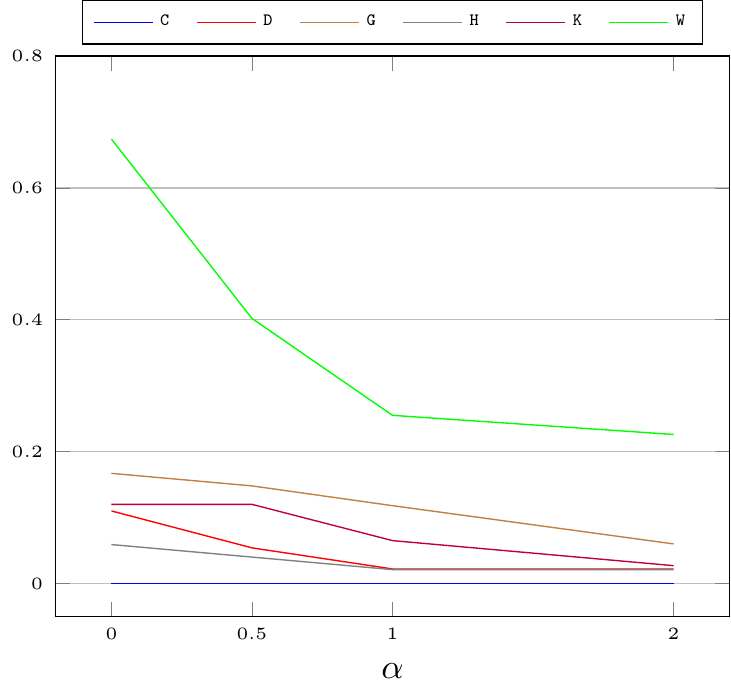}
		\caption{Price of efficiency (PoE).}
	\end{subfigure}
	\caption{PoF and PoE averaged for each value of $\alpha$ for \texttt{Continuous} location problem.}
	\label{f:continuous_np}
\end{figure}

Figures \ref{f:discrete_np} and \ref{f:continuous_np} show the values of the price of fairness (PoF) and the price of efficiency (PoE) defined in Section \ref{sec:FOWA} for the discrete and continuous location problems, respectively.
These values are shown for each type of problem ($\lambda$) considered and for each value of $\alpha$. 
The highest values of the PoF are obtained with \texttt{C} (the one for which PoE is the smallest) whereas the highest values of the PoE are obtained with \texttt{W} (the one with smallest PoF). However, when the parameter $\alpha$ increases, we obtain a range of fair solutions between the classical MCLP ($\alpha = 0,\; \lambda = \texttt{W}$) and the maximin approach (\texttt{C}), balancing efficiency and fairness.

The performance of the solutions for each of the problem types $\lambda$ seems to indicate a relationship between the PoF and the PoE, in the sense that the largest PoF is related to the smallest PoE (or vice versa) as can be seen for \texttt{C} (or \texttt{W}).
Nevertheless, this is not true in general. For instance, in the discrete case  (Figure \ref{f:discrete_np}), we got that \texttt{D} have a value of PoF between those obtained by \texttt{G} and \texttt{K} but its PoE is below \texttt{K}. That is, the solution for \texttt{D} is better than  \texttt{K} both in efficiency and fairness.

Finally, for the discrete location problem we can conclude that the trade-off between fairness and efficiency is promising when considering as $\lambda \in \{\texttt{D},\texttt{G}\}$, obtaining fairer solutions by increasing the value of $\alpha$. If one is willing to lose more efficiency to gain fairness the better operator seems to be the harmonic (\texttt{H}). 

For the continuous problems, similar conclusions are derived.  Specifically, the harmonic operator (\texttt{H}) seems to provide fair solutions with values of PoE close to zero, but still with not so high values of PoF as \texttt{C}.

These graphs show that for different values the decision maker has a wide range of possibilities to decide in choosing between the fairest or the one that covers the most.
This is the first conclusion that we can draw from the results and that verifies the starting hypothesis, a new fairness measurement function for allocation problems.

More detailed graphs of the averaged results can be found in  \ref{sec:appendix}.

One can draw similar conclusions to those presented above. Nevertheless, one can draw other interesting observations from the averaged results.  For instance, in Figure \ref{f:d-PoF-n}, for $n=45$, a breakpoint can be observed for \texttt{H}  in $\alpha = 1$. Observe that for small values of $n$ and large value of $p$, only a few instances were solved with MIPGap less than 5\%, and then, the  maximum averaged value of PoF is reached.

Another behavior that needs further explanation is for when we average by $p$ for both type of location spaces, and for both PoF and PoE. For $p = 20$, we observe that for $\lambda \in \{\texttt{H}, \texttt{K}\}$ the graphs drop to 0 for $\alpha \in \{1,2\}$ in discrete case and for $\lambda = \texttt{G}$ in continuous case. In this case, none of the instances were solved with MIPGap under $\%5$. Moreover, for the continuous case with $p = 20$, Figure \ref{f:c-PoF-p} shows the PoF obtained for \texttt{H} is 0 for all $\alpha$ values, i.e., it fails to improve the classical MCLP solution, but the PoE shows in Figure \ref{f:c-PoE-p} that when $\alpha \in \{0,0.5\}$ a better solution in terms of efficiency was obtained, even if the coverage is equal.  This  empirically proof that one can find alternative maximal coverage solutions with different allocation of the users of customers, one of them fairer than others.

\begin{figure}[H]
	\begin{subfigure}[b]{.45\linewidth}
		\centering 
		\includegraphics[scale=0.95]{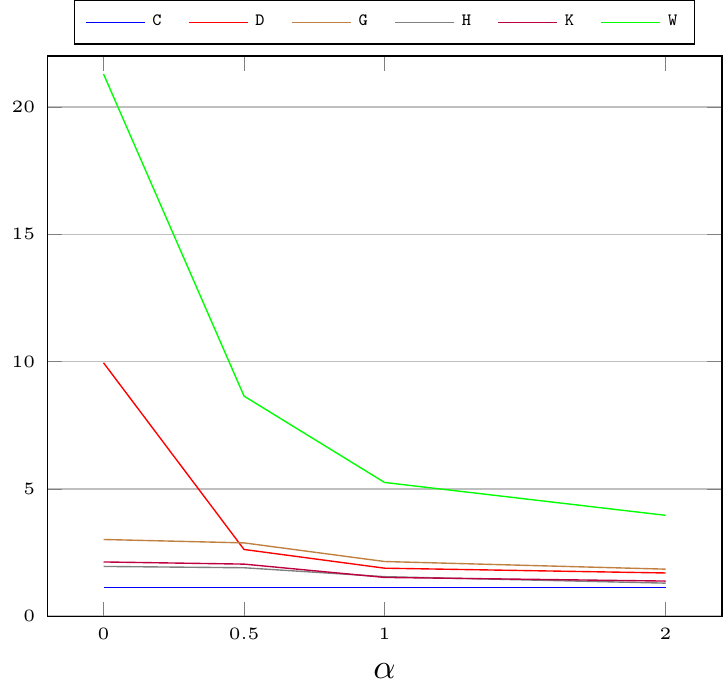}
		\caption{\texttt{Disc} location problem.}
	\end{subfigure}~\begin{subfigure}[b]{.45\linewidth}
		\centering 
		\includegraphics[scale=0.95]{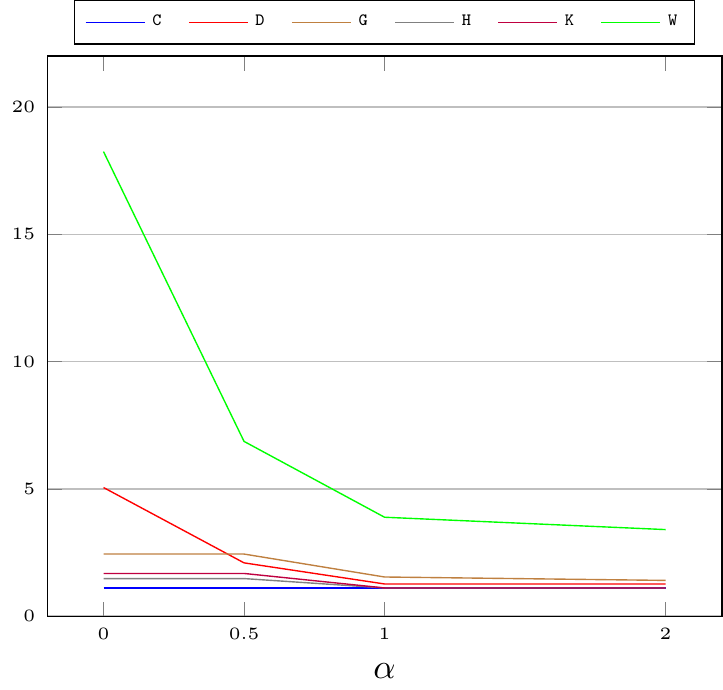}
		\caption{\texttt{Cont} location problem.}
	\end{subfigure}
	\caption{Gini averaged for each value of $\alpha$ for both types of location spaces.}
	\label{f:gini}
\end{figure}

We have also analyzed the behavior of the solutions obtained with our approaches with respect to the Gini index (Section \ref{sec:FOWA}). The Gini index measures the global pairwise envy of the different facilities (in this case the  larger the difference in covered demand between two facilities, the larger the Gini index). In Figure \ref{f:gini} we show the overall performance of all the solved instances in terms of the average Gini indices for each of the used $\lambda$-weights and each of the values of $\alpha$. Observe that the behavior is very similar to PoE. Again, the \texttt{C} models outperforms the rest of models for all $\alpha$, while \texttt{W} exhibits the worse behavior, being the maximum-envy solution. In the middle, \texttt{D} provides a trade off and \texttt{H} has also a better performance than \texttt{W}.  For the continuous location case, both \texttt{K} and \texttt{H} have similar Gini indices than \texttt{C}, in contrast to the PoE. This is because the difference between the facilities is the same as \texttt{C} and therefore minimal, but does not reach the minimum coverage.

\begin{figure}[H]
	\begin{subfigure}[b]{.32\linewidth}
		\centering \fbox{\parbox[c][43mm][b]{40mm}{\centering \includegraphics[scale = 0.19]{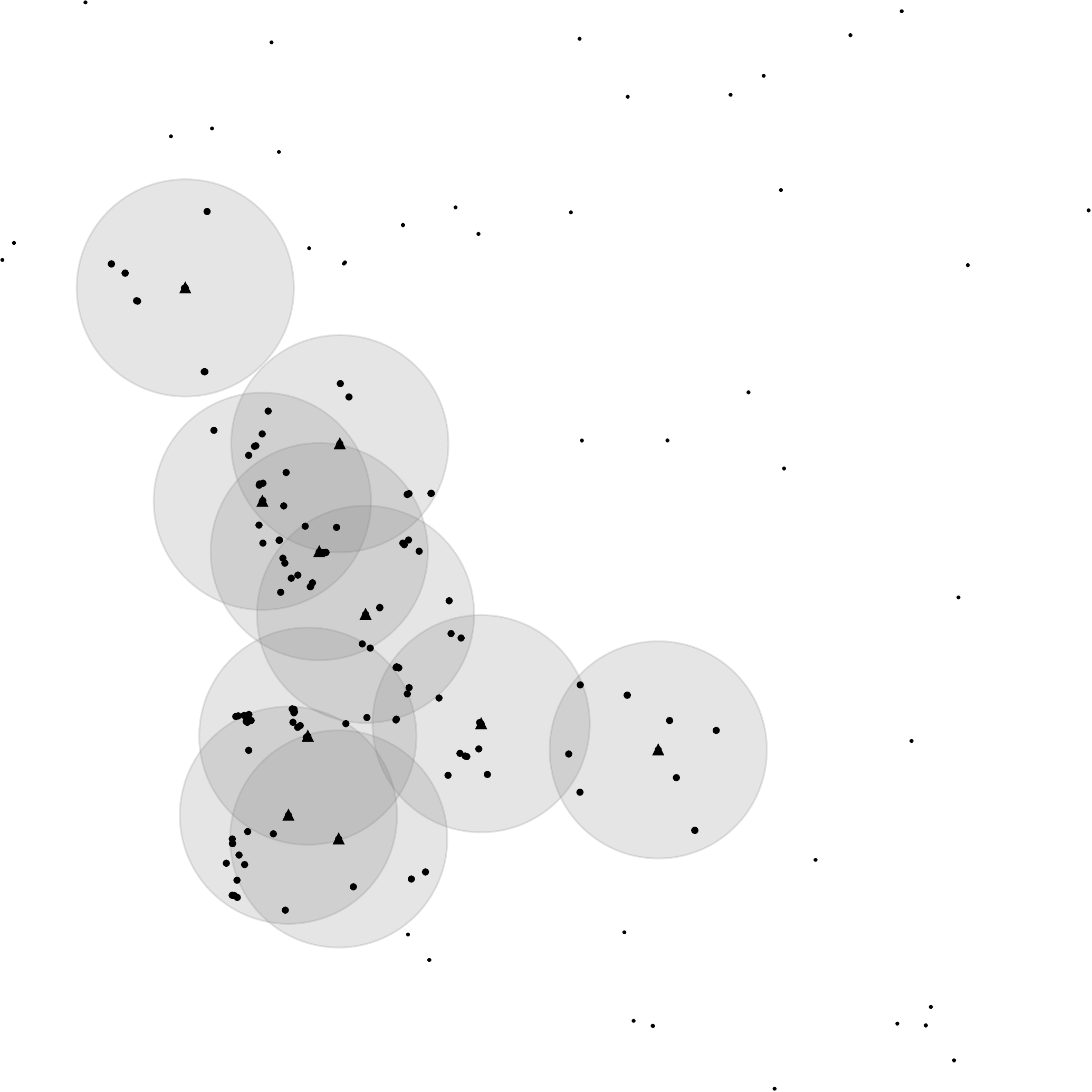}}}
		\caption{\texttt{C} (62.46\%)}\label{figa}
	\end{subfigure}\hfill 
	\begin{subfigure}[b]{.32\linewidth}
		\centering\fbox{\parbox[c][43mm][b]{40mm}{\centering \includegraphics[scale = 0.19]{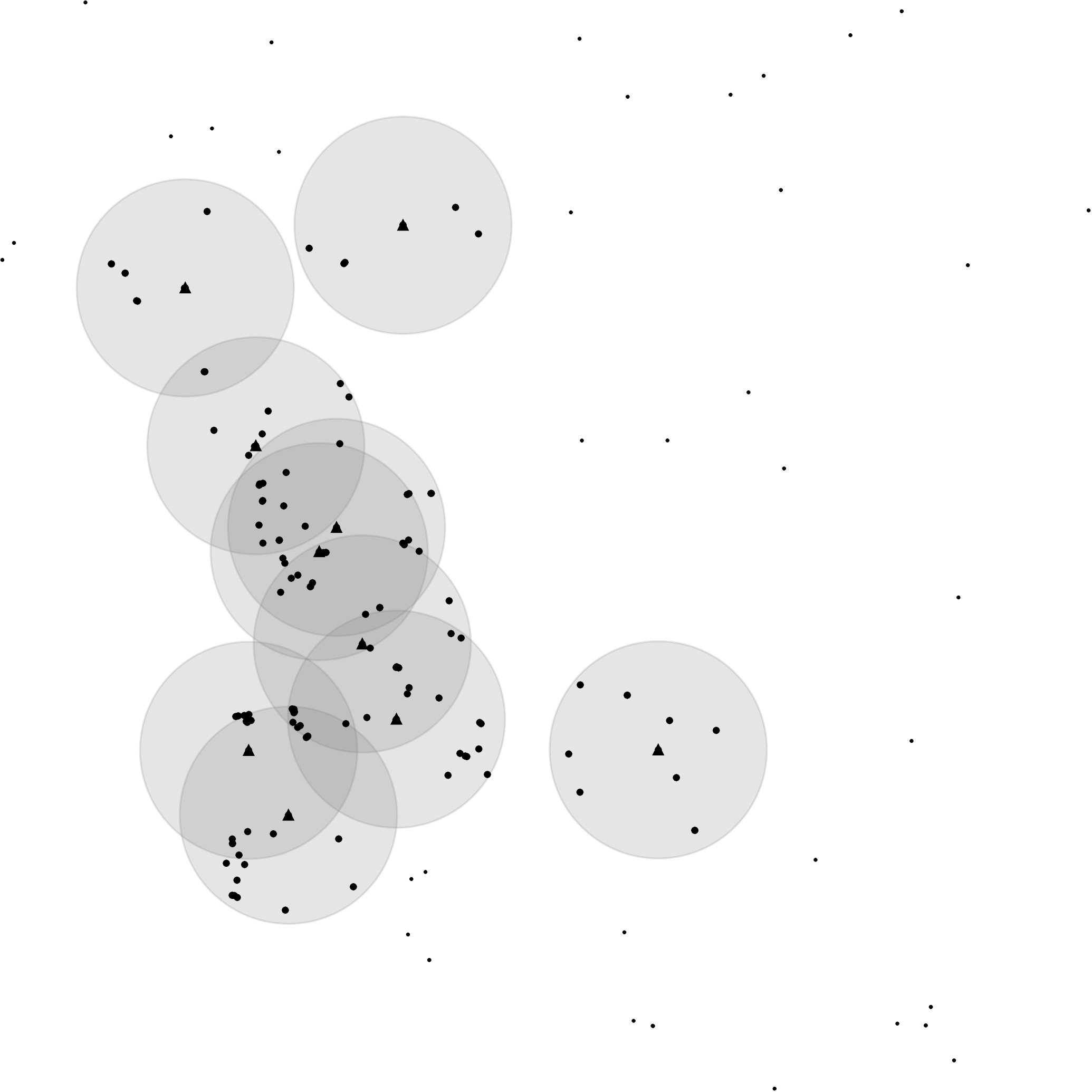}}}
		\caption{\texttt{K} and $\alpha = 1$ (67.18\%)}\label{figb}
	\end{subfigure}\hfill 
	\begin{subfigure}[b]{.32\linewidth}
		\centering\fbox{\parbox[c][43mm][b]{40mm}{\centering \includegraphics[scale = 0.19]{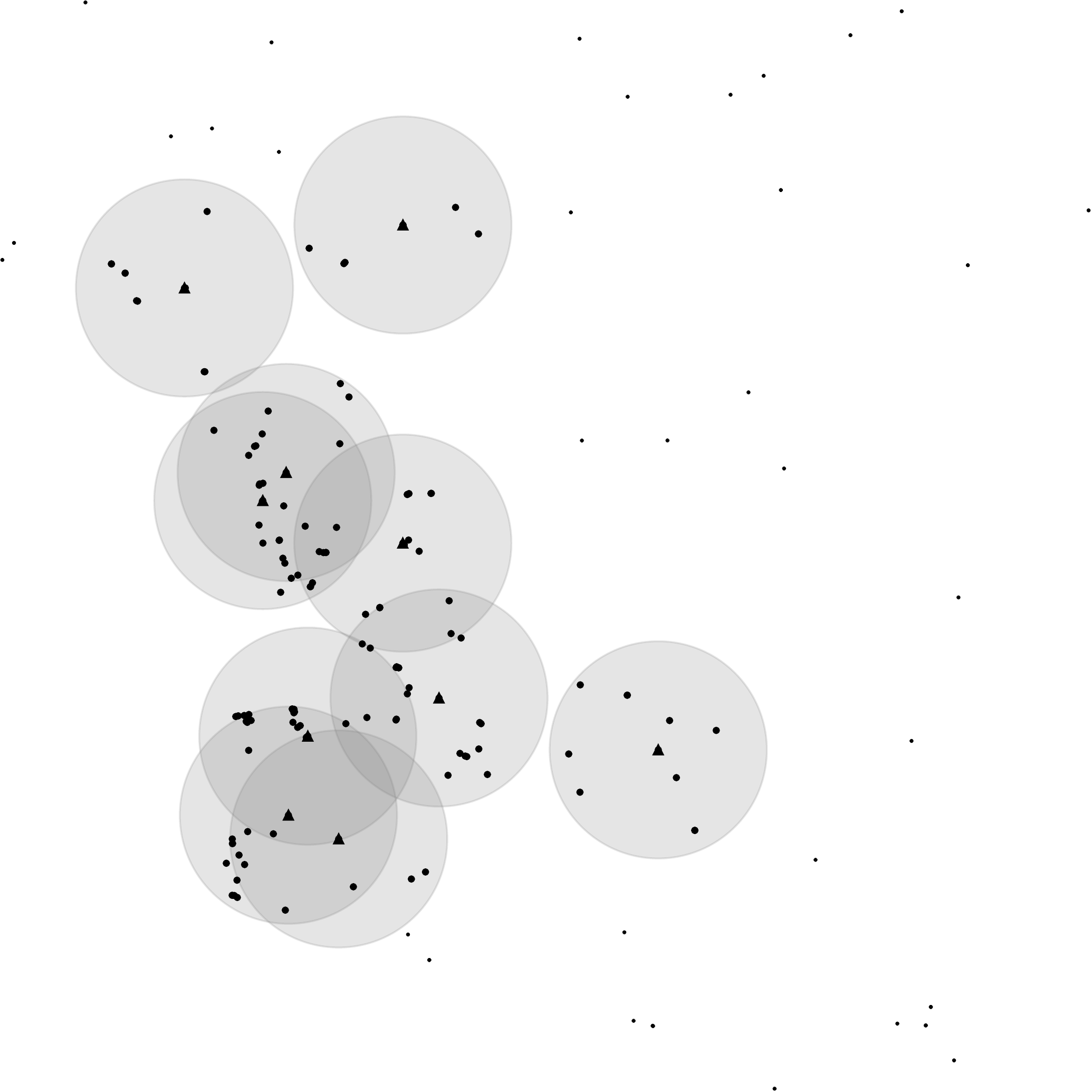}}}
		\caption{\texttt{D} and $\alpha = 0.5$ (67.79\%)}\label{figc}
	\end{subfigure}\\ \quad \\
\begin{subfigure}[b]{.32\linewidth}
	\centering \fbox{\parbox[c][43mm][b]{40mm}{\centering \includegraphics[scale = 0.19]{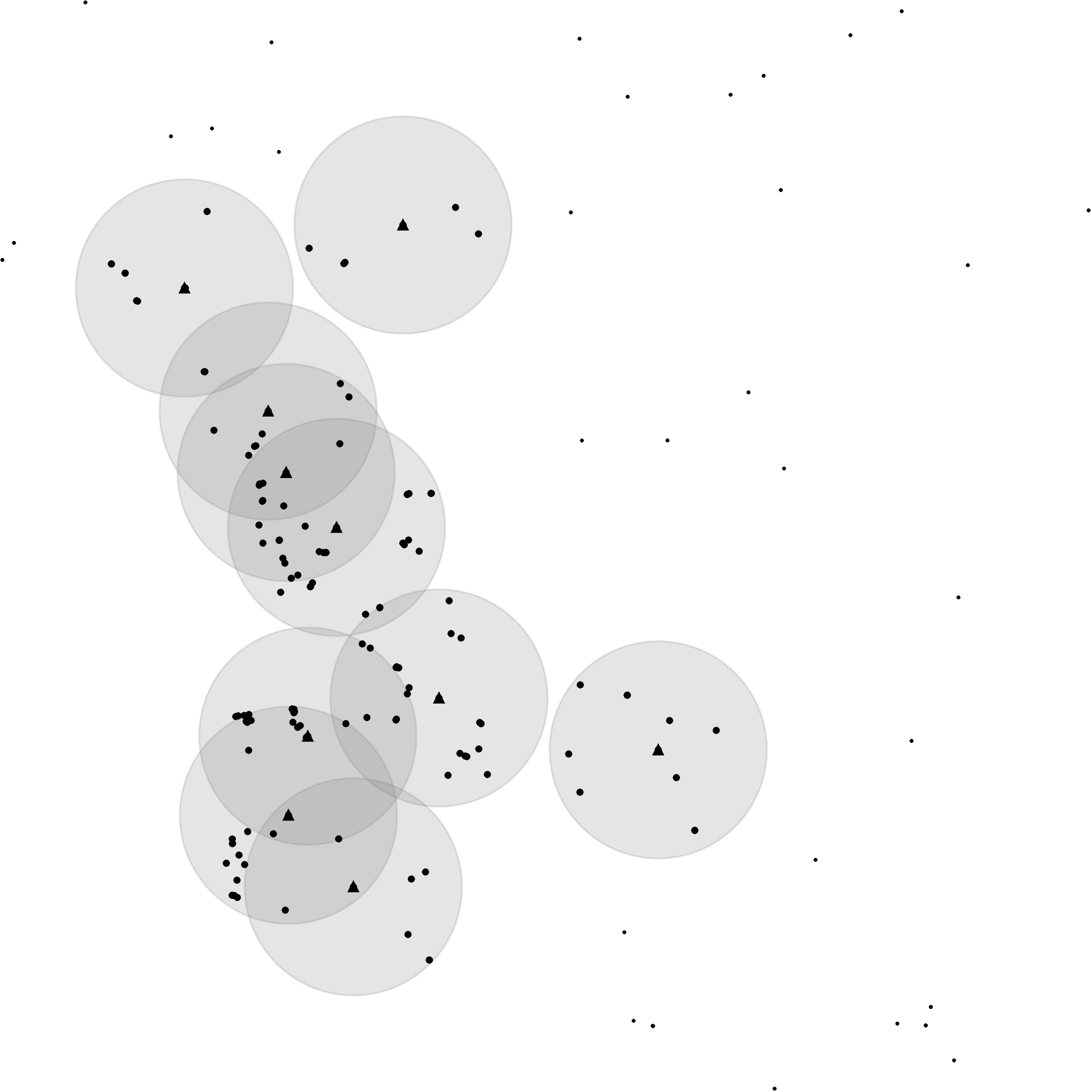}}}
	\caption{\texttt{D} and $\alpha = 0$ (70.02\%)}\label{figd}
\end{subfigure}\hfill 
\begin{subfigure}[b]{.32\linewidth}
	\centering\fbox{\parbox[c][43mm][b]{40mm}{\centering \includegraphics[scale = 0.19]{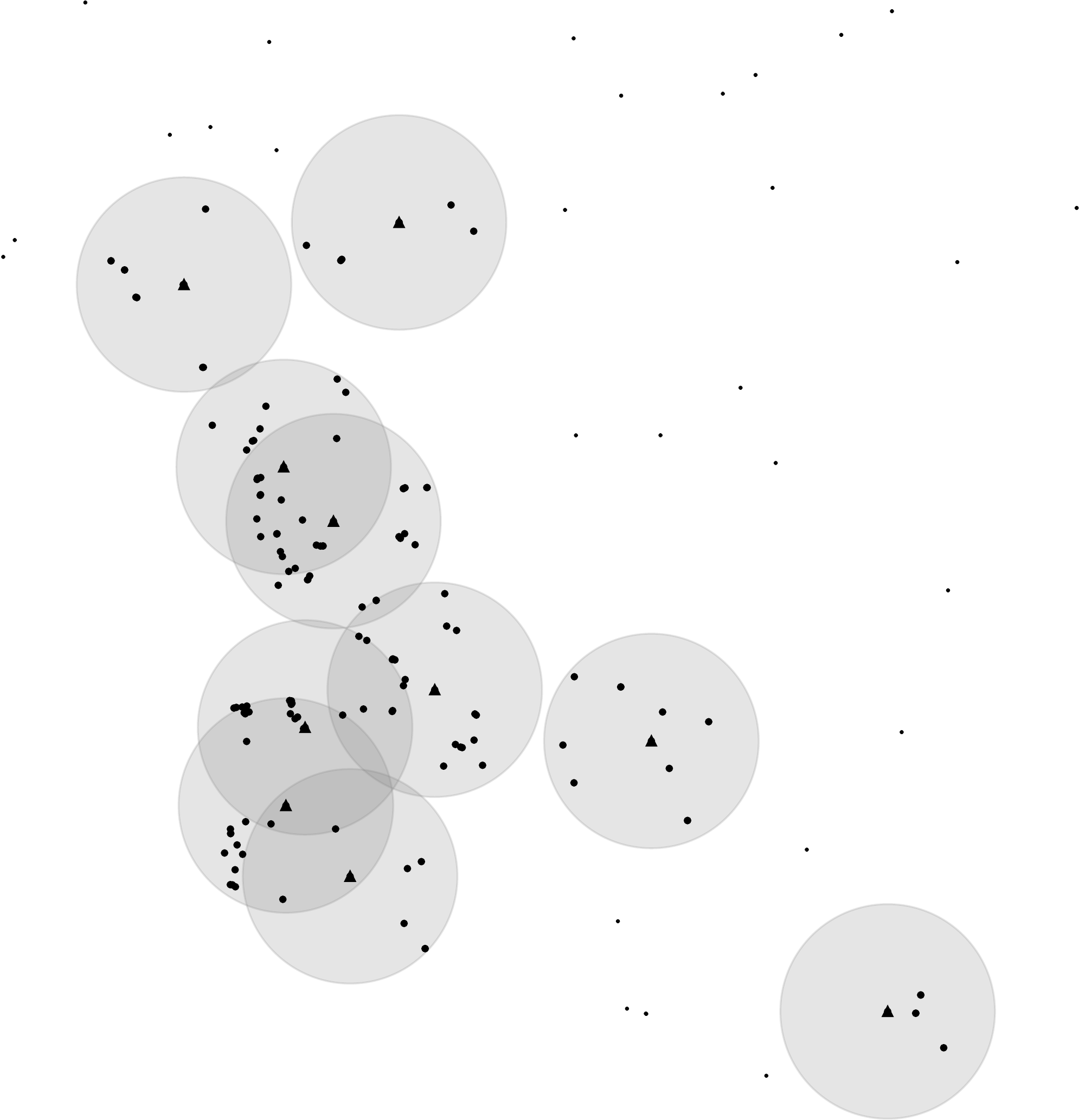}}}
	\caption{\texttt{W} and $\alpha = 0.5$ (73.13\%)}\label{fige}
\end{subfigure}\hfill 
\begin{subfigure}[b]{.32\linewidth}
	\centering\fbox{\parbox[c][43mm][b]{40mm}{\centering \includegraphics[scale = 0.19]{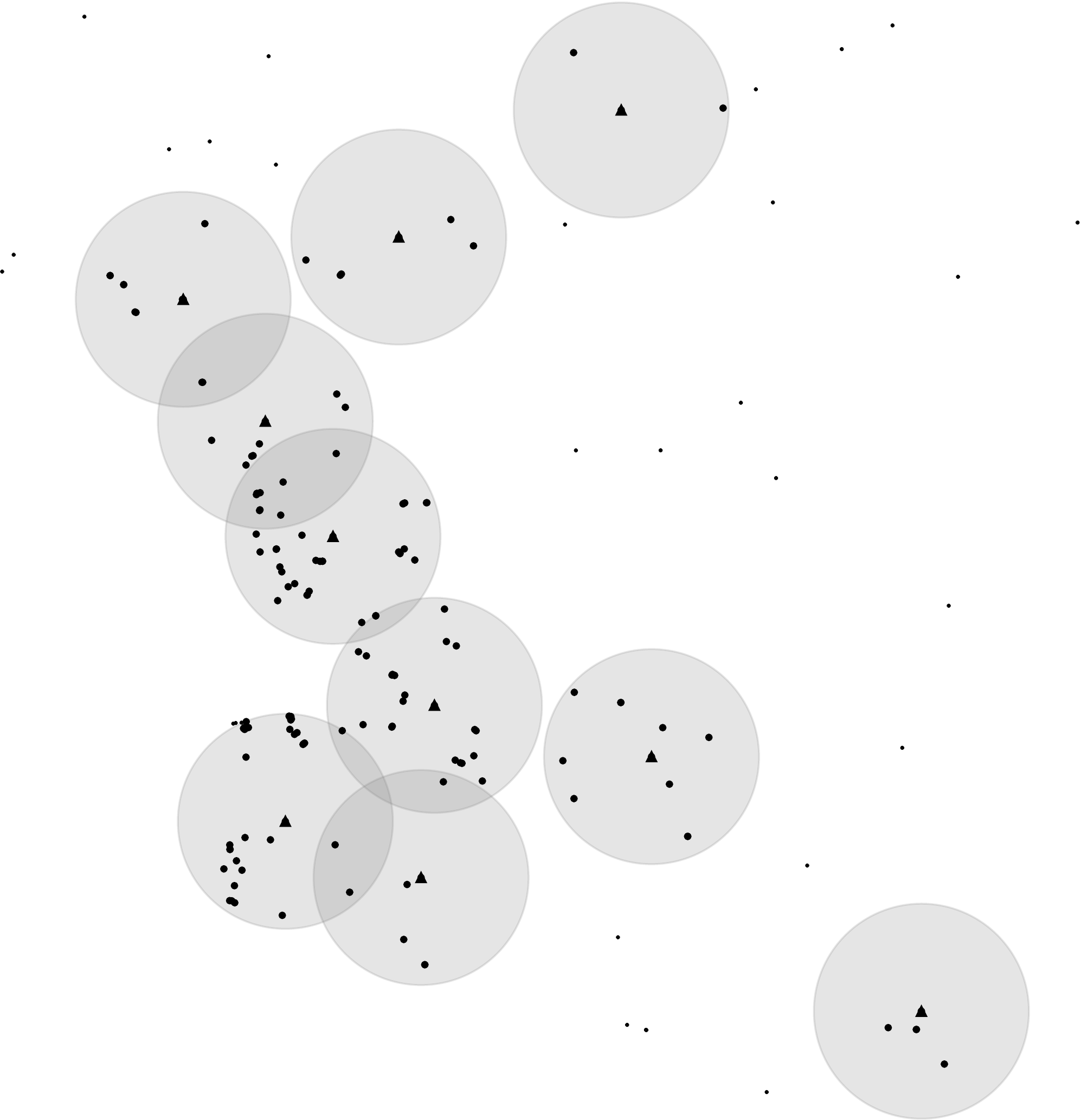}}}
	\caption{\texttt{W} and $\alpha = 0$ (74.09\%)}\label{figf}
\end{subfigure}
	\caption{Example of fairness distribution for some results when $n=179$ and $p=10$, and in parenthesis the weighted covered percentage.}\label{f:facilities}
\end{figure}

Finally, in Figure \ref{f:facilities} we illustrate different optimal solutions obtained for a selected instance for different values of ($\lambda$, $\alpha$) of fairness measure.  From left top to right bottom we draw:  (\texttt{C}, any $\alpha$), $(\texttt{K},1)$,  $(\texttt{D},0.5)$, $(\texttt{D},0)$, $(\texttt{W},0.5)$ and $(\texttt{W},0)$. The reader can observe that, geometrically, a greater concentration of facilities is obtained for fairer solutions, being not only the position of the facilities the factor that affect the fairness of a solution, but also the allocation of the users to them, which can be more easily when different facilities are able to cover the same demand points.
	
%%%% Figuras tesis

%\begin{figure}[H]
%	\begin{subfigure}[b]{.32\linewidth}
%		\centering \fbox{\parbox[c][43mm][b]{36mm}{\centering \includegraphics[scale = 0.2]{C_0.00.pdf}}}
%		\caption{\texttt{C}}\label{figa}
%	\end{subfigure}\hfill 
%	\begin{subfigure}[b]{.32\linewidth}
%		\centering\fbox{\parbox[c][43mm][b]{36mm}{\centering \includegraphics[scale = 0.2]{H_0.00.pdf}}}
%		\caption{\texttt{H} and $\alpha = 0$}\label{figb}
%	\end{subfigure}\hfill 
%	\begin{subfigure}[b]{.32\linewidth}
%		\centering\fbox{\parbox[c][43mm][b]{36mm}{\centering \includegraphics[scale = 0.2]{D_0.50.pdf}}}
%		\caption{\texttt{D} and $\alpha = 0.5$}\label{figc}
%	\end{subfigure}\\ \quad \\
%	\begin{subfigure}[b]{.32\linewidth}
%		\centering \fbox{\parbox[c][43mm][b]{36mm}{\centering \includegraphics[scale = 0.2]{G_0.00.pdf}}}
%		\caption{\texttt{G} and $\alpha = 0$}\label{figd}
%	\end{subfigure}\hfill 
%	\begin{subfigure}[b]{.32\linewidth}
%		\centering\fbox{\parbox[c][43mm][b]{36mm}{\centering \includegraphics[scale = 0.2]{W_0.50.pdf}}}
%		\caption{\texttt{W} and $\alpha = 0.5$}\label{fige}
%	\end{subfigure}\hfill 
%	\begin{subfigure}[b]{.32\linewidth}
%		\centering\fbox{\parbox[c][43mm][b]{36mm}{\centering \includegraphics[scale = 0.2]{W_0.00.pdf}}}
%		\caption{\texttt{W} and $\alpha = 0$}\label{figf}
%	\end{subfigure}
%	\caption{Example of fairness distribution for some results when $n=90$ and $p=10$.}\label{f:facilities}
%\end{figure}

\section{Conclusions and further research}

We present in this paper a novel fairness measure for Maximal Covering Location Problems, that combines the OWA operators, early introduced by \cite{yager88}, and the $\alpha$-fairness scheme introduced by \cite{atkinson70}. 

We develop suitable mathematical optimization models that allow to capture the notion of fairness in the MCLP for the two main types of location spaces that are studied in the literature: the discrete and continuous. The models are then reformulated as MISOCO problem, using the geometrical insights of the problem, and then, the programs that we propose are able to be solved with the available off-the-shelf software.

We have tested our models using a real data set containing the locations of residential schools and student hostels in Canada. Applying our new scheme to this dataset, we empirically observed that our models provide different solutions in terms of fairness than using the OWA operators and the $\alpha$-fairness scheme separately. Therefore, we conclude that the combination of both in the same allocation scheme provides the decision maker with a wide range of options to find a trade-off between fairness and efficiency.

Further research lines on fairness topic include the using of this combination of schemes into different location problems, as the set covering problem, location problems with capacities and also its incorporation to queue problems where the facilities has to deal with the management of waiting users when the facility is saturated, and fair solutions of the problem could be a successfully too to deal with this problematic.

\section*{Acknowledgments}

This research has been partially supported by Spanish Ministerio de Ciencia e Innovación, AEI/FEDER grant number PID2020-114594GBC21, Junta de Andalucía projects P18-FR-1422/2369 and projects FEDERUS-1256951, B-FQM-322-UGR20, CEI-3-FQM331 and NetmeetData (Fundación BBVA 2019). The first author was also partially supported by the IMAG-Maria de Maeztu grant CEX2020-001105-M /AEI /10.13039/501100011033.
The second author was also partially supported by Spanish Ministry of Education and Science grant number PEJ2018-002962-A, and the PhD Program in Mathematics at the Universidad de Granada.

\appendix

\section{Discussions}\label{sec:appendix}

\begin{figure}[H]
	\begin{center}
		\begin{subfigure}[b]{.45\linewidth}
			\includegraphics[scale=1]{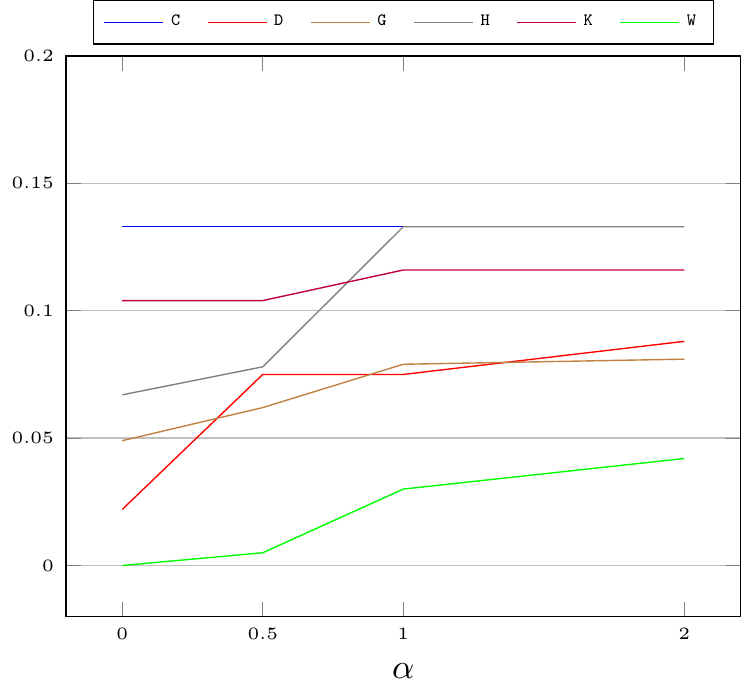}	
			\caption{$n=45$.}
		\end{subfigure}~
		\begin{subfigure}[b]{.45\linewidth}
			\includegraphics[scale=1]{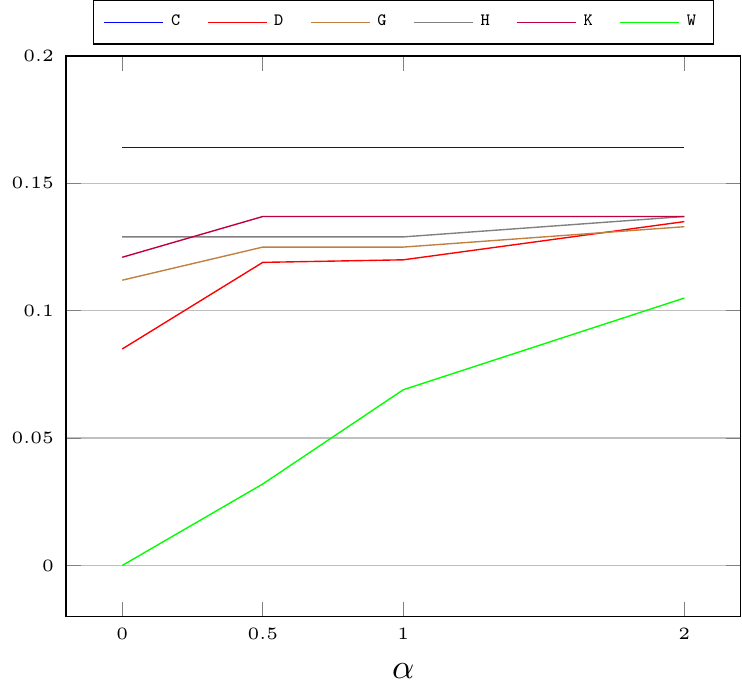}	
			\caption{$n=90$.}
		\end{subfigure}
		
		\begin{subfigure}[b]{.45\linewidth}
			\includegraphics[scale=1]{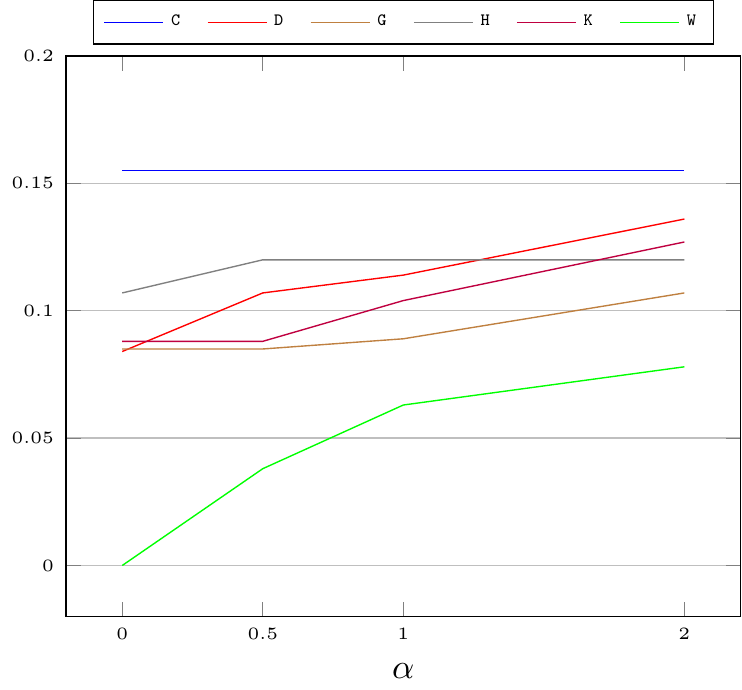}	
			\caption{$n=120$.}
		\end{subfigure}~
		\begin{subfigure}[b]{.45\linewidth}
			\includegraphics[scale=1]{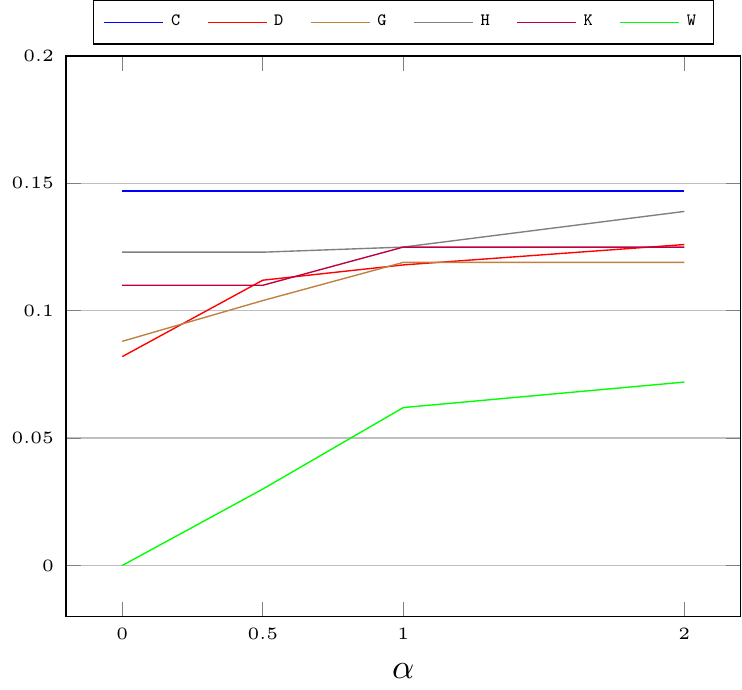}	
			\caption{$n=179$.}
		\end{subfigure}~
	\end{center}
	\caption{Price of fairness averaged by $\alpha \in \{0,0.5,1,2\}$ and $n \in\{45,90,120,179\}$ \texttt{Discrete} location problem.}
	\label{f:d-PoF-n}
\end{figure}

\begin{figure}[H]
	\begin{center}
		\begin{subfigure}[b]{.45\linewidth}
			\includegraphics[scale=1]{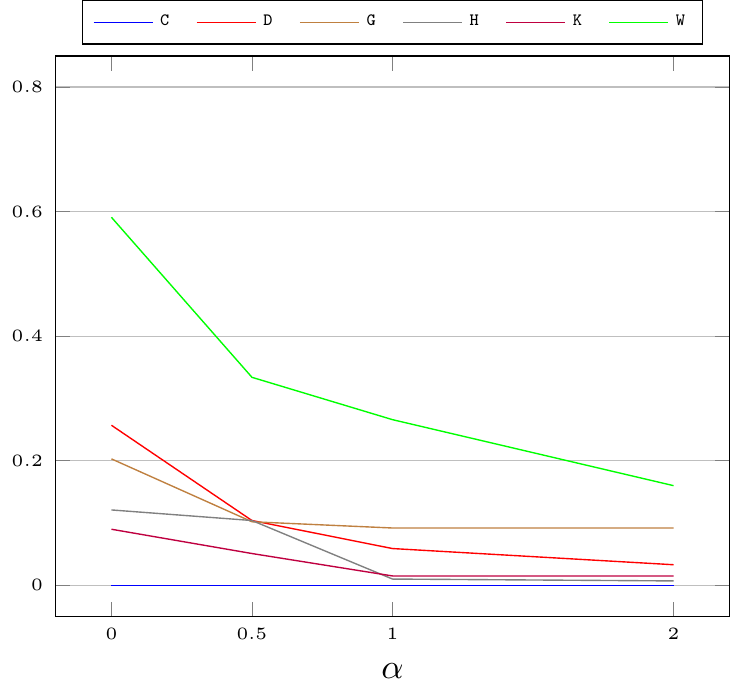}	
			\caption{$n=45$.}
		\end{subfigure}~
		\begin{subfigure}[b]{.45\linewidth}
			\includegraphics[scale=1]{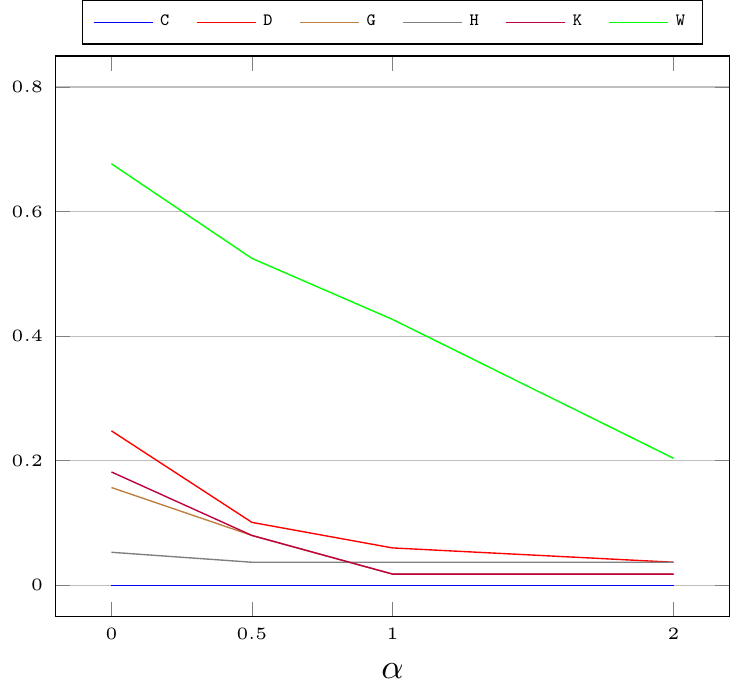}	
			\caption{$n=90$.}
		\end{subfigure}
		
		\begin{subfigure}[b]{.45\linewidth}
			\includegraphics[scale=1]{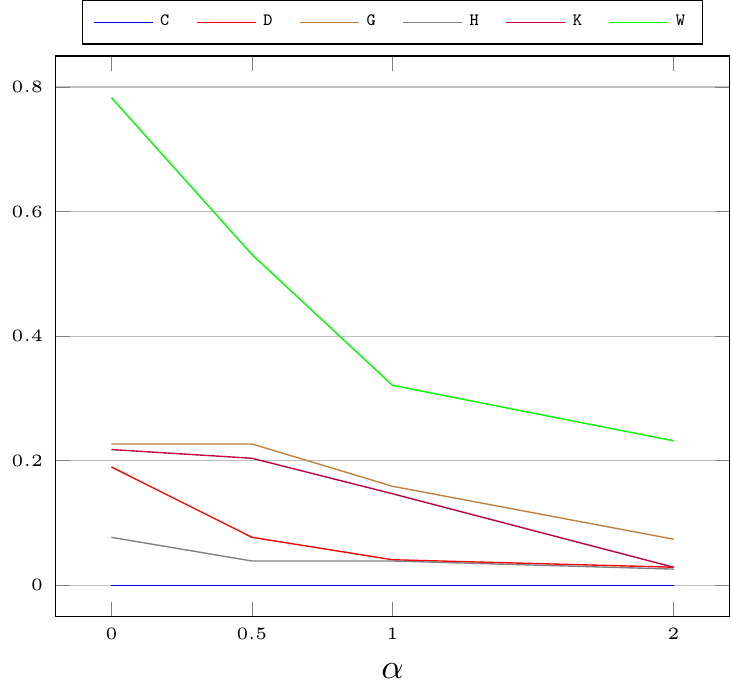}	
			\caption{$n=120$.}
		\end{subfigure}~
		\begin{subfigure}[b]{.45\linewidth}
			\includegraphics[scale=1]{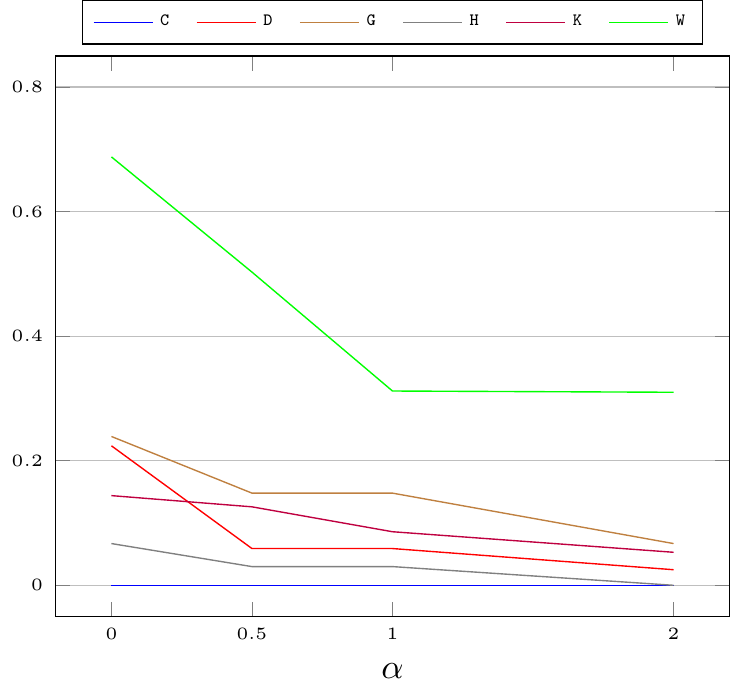}	
			\caption{$n=179$.}
		\end{subfigure}~
	\end{center}
	\caption{Price of efficiency averaged by $\alpha \in \{0,0.5,1,2\}$ and $n \in\{45,90,120,179\}$ for \texttt{Discrete} location problem.}
	\label{f:d-PoE-n}
\end{figure}

\begin{figure}[H]
	\begin{center}
				\begin{subfigure}[b]{.45\linewidth}
			\includegraphics[scale=1]{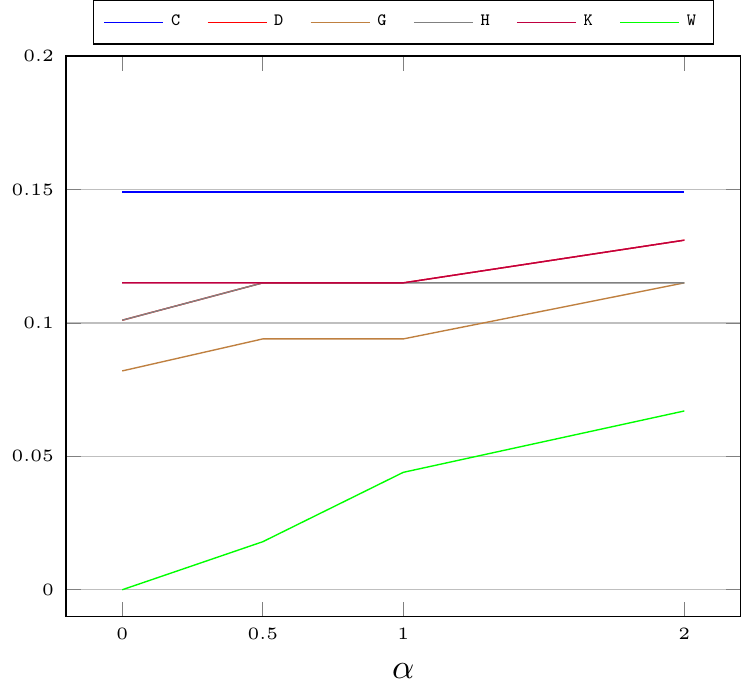}	
			\caption{$p=5$.}
		\end{subfigure}~
		\begin{subfigure}[b]{.45\linewidth}
			\includegraphics[scale=1]{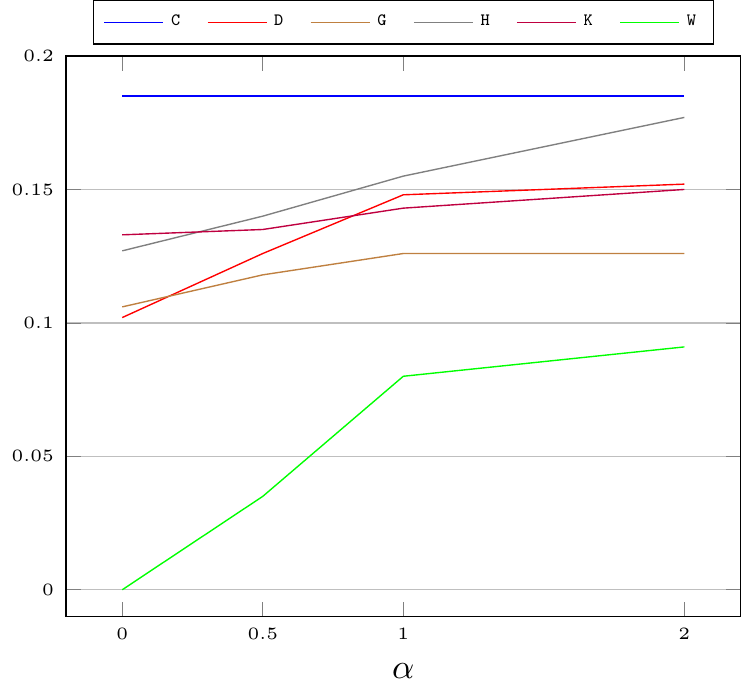}	
			\caption{$p=10$.}
		\end{subfigure}
		
		\begin{subfigure}[b]{.45\linewidth}
			\includegraphics[scale=1]{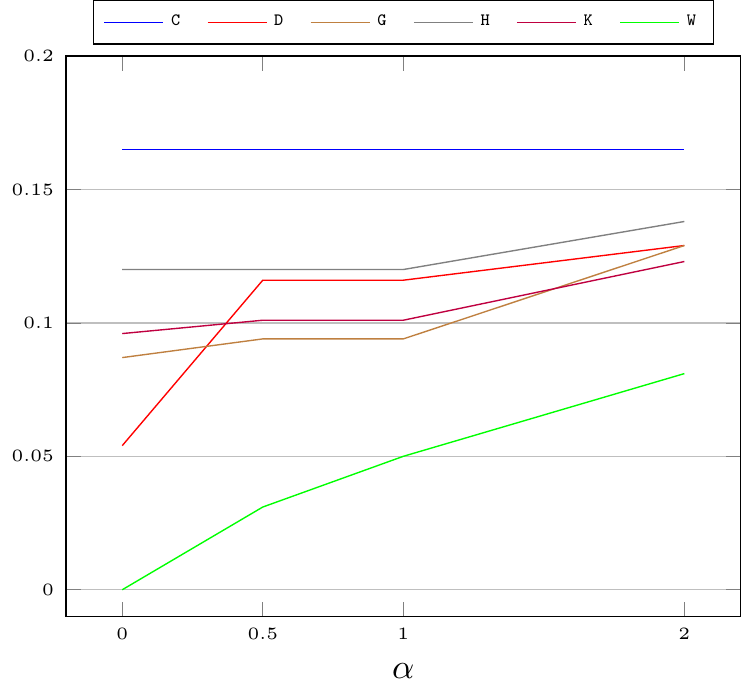}	
			\caption{$p=15$.}
		\end{subfigure}~
		\begin{subfigure}[b]{.45\linewidth}
			\includegraphics[scale=1]{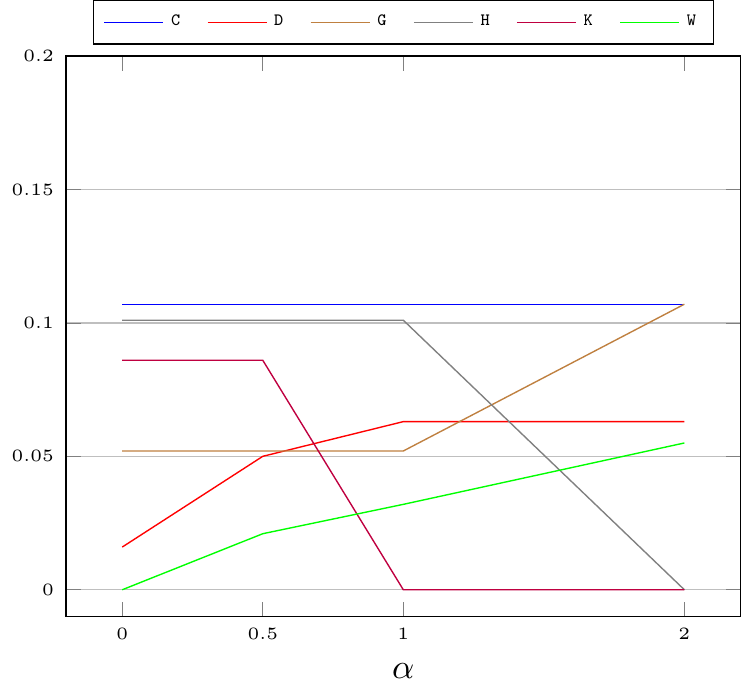}	
			\caption{$p=20$.}
		\end{subfigure}
	\end{center}
	\caption{Price of fairness averaged by $\alpha \in \{0,0.5,1,2\}$ and $p \in\{5,10,15,20\}$ \texttt{Discrete} location problem.}
	\label{f:d-PoF-p}
\end{figure}

\begin{figure}[H]
	\begin{center}
		\begin{subfigure}[b]{.45\linewidth}
			\includegraphics[scale=1]{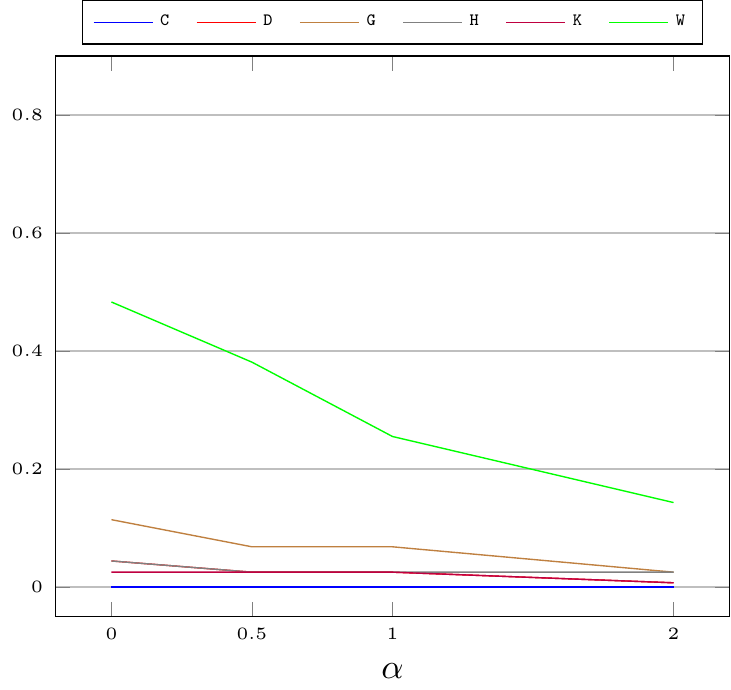}	
			\caption{$p=5$.}
		\end{subfigure}~
		\begin{subfigure}[b]{.45\linewidth}
			\includegraphics[scale=1]{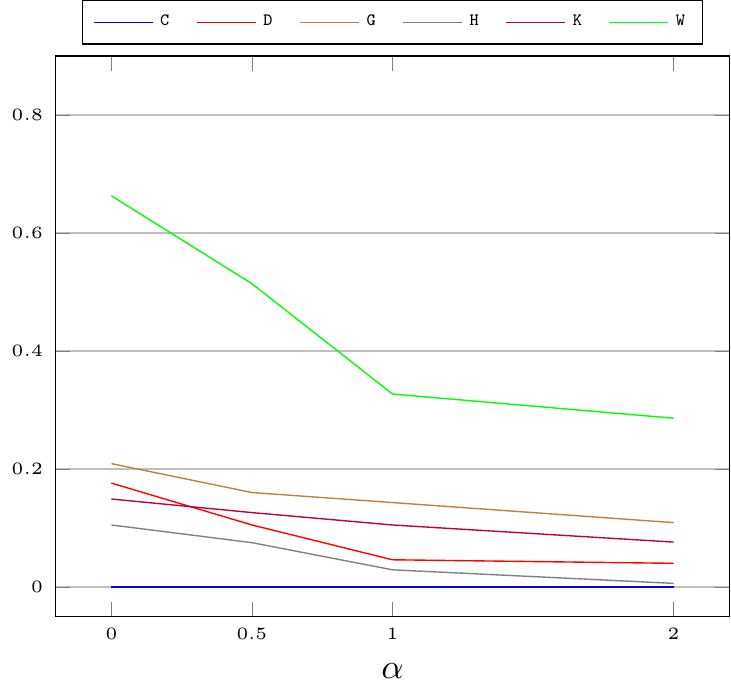}	
			\caption{$p=10$.}
		\end{subfigure}
		
		\begin{subfigure}[b]{.45\linewidth}
			\includegraphics[scale=1]{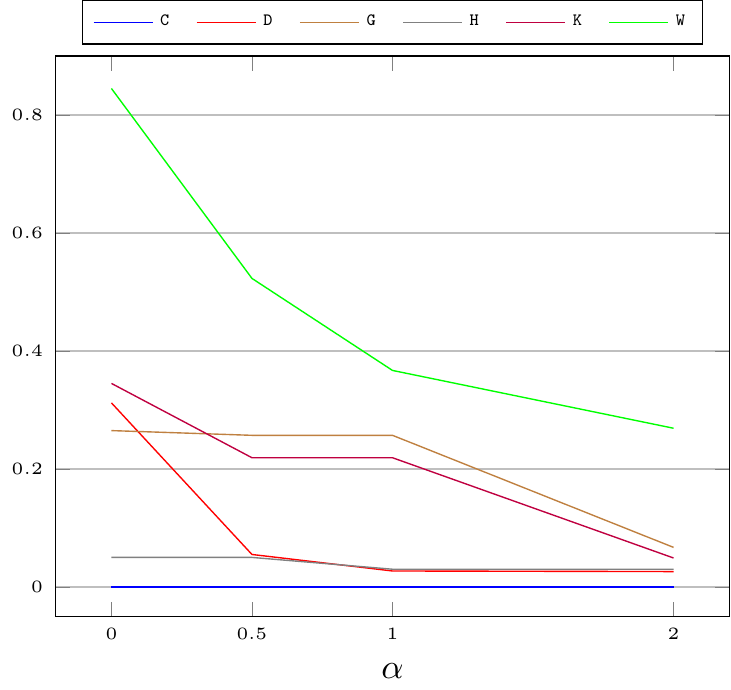}	
			\caption{$p=15$.}
		\end{subfigure}~
		\begin{subfigure}[b]{.45\linewidth}
			\includegraphics[scale=1]{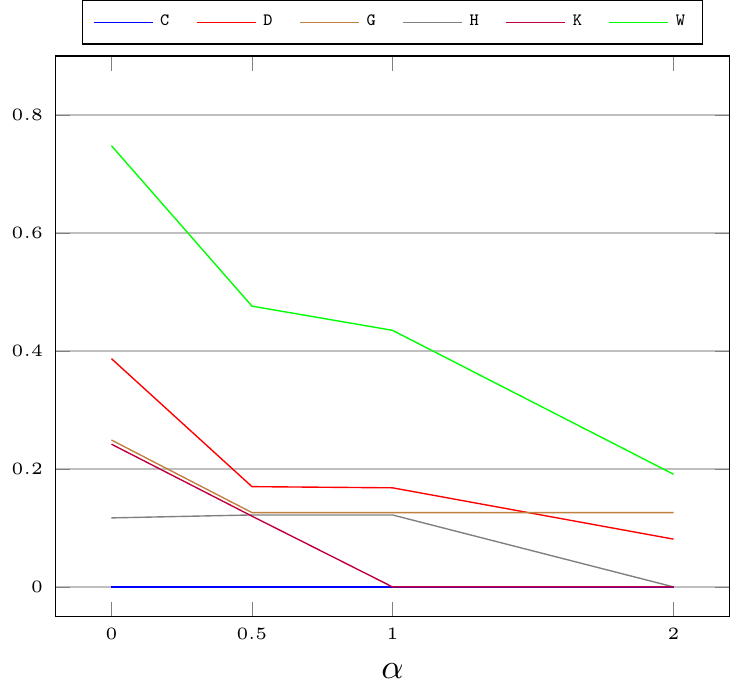}	
			\caption{$p=20$.}
		\end{subfigure}
	\end{center}
	\caption{Price of efficiency averaged by $\alpha \in \{0,0.5,1,2\}$ and $p \in\{5,10,15,20\}$ for \texttt{Discrete} location problem.}
	\label{f:d-PoE-p}
\end{figure}

\begin{figure}[H]
	\begin{center}
		\begin{subfigure}[b]{.45\linewidth}
			\includegraphics[scale=1]{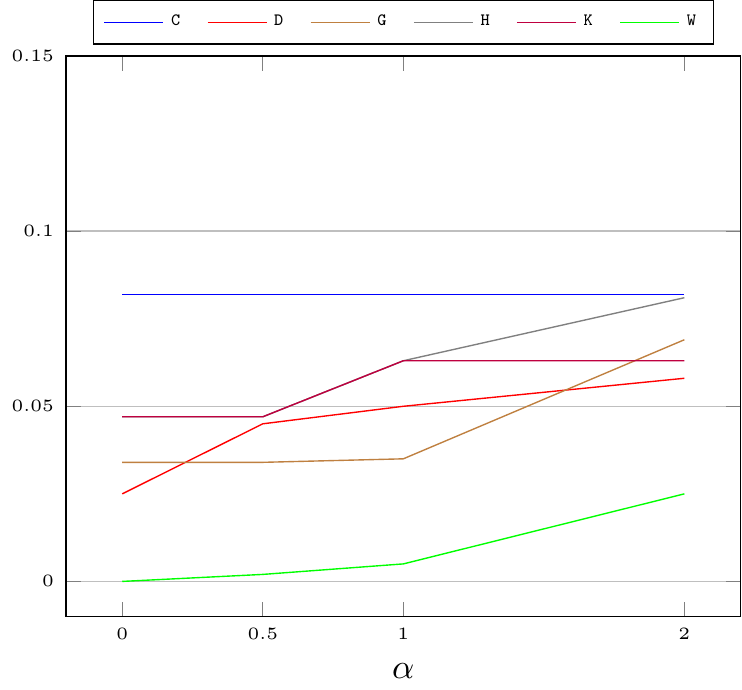}	
			\caption{$n=45$.}
		\end{subfigure}~
		\begin{subfigure}[b]{.45\linewidth}
			\includegraphics[scale=1]{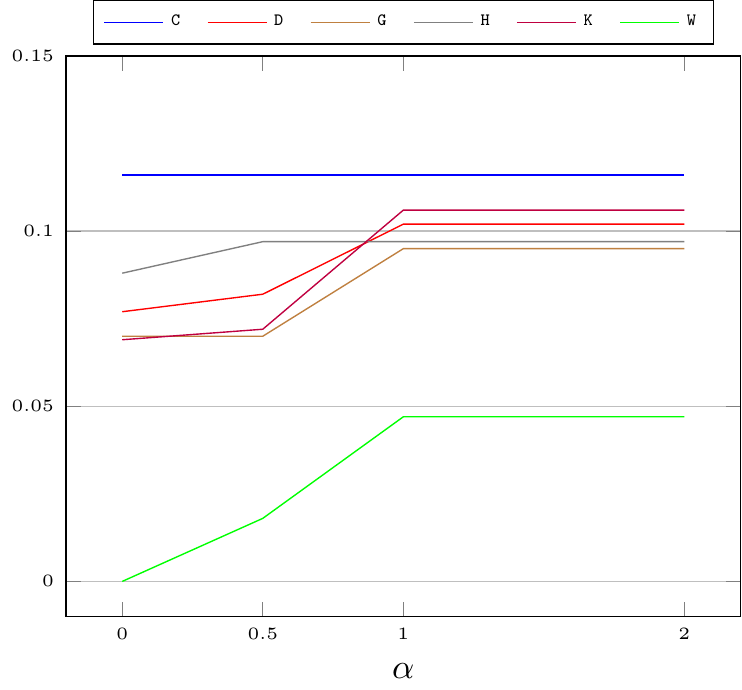}	
			\caption{$n=90$.}
		\end{subfigure}
		
		\begin{subfigure}[b]{.45\linewidth}
			\includegraphics[scale=1]{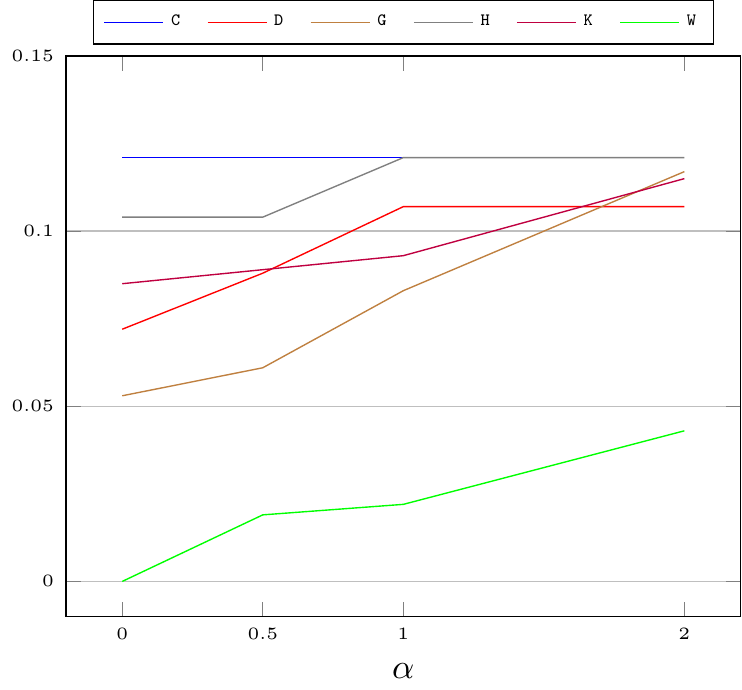}	
			\caption{$n=120$.}
		\end{subfigure}~
		\begin{subfigure}[b]{.45\linewidth}
			\includegraphics[scale=1]{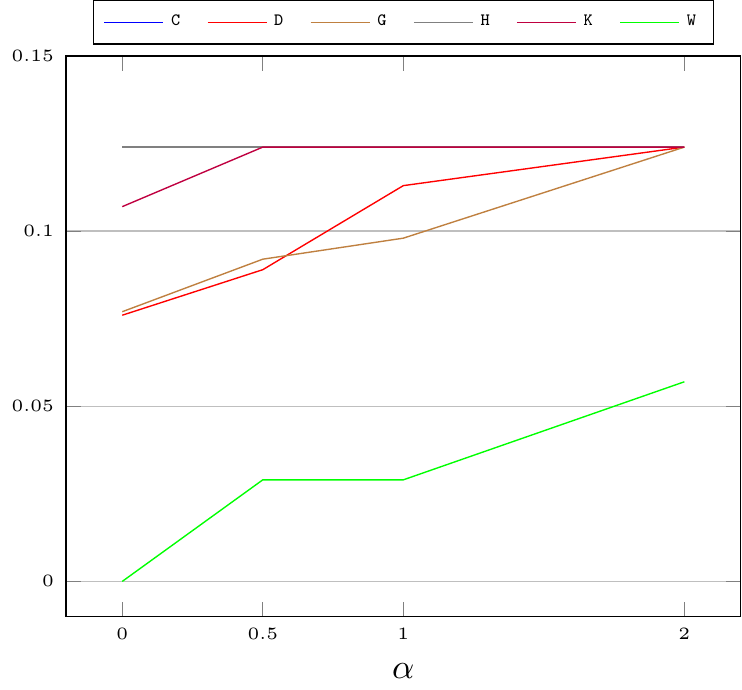}	
			\caption{$n=179$.}
		\end{subfigure}~
	\end{center}
	\caption{Price of fairness averaged by $\alpha \in \{0,0.5,1,2\}$ and $n \in\{45,90,120,179\}$ \texttt{Continuous} location problem.}
	\label{f:c-PoF-n}
\end{figure}

\begin{figure}[H]
	\begin{center}
		\begin{subfigure}[b]{.45\linewidth}
			\includegraphics[scale=1]{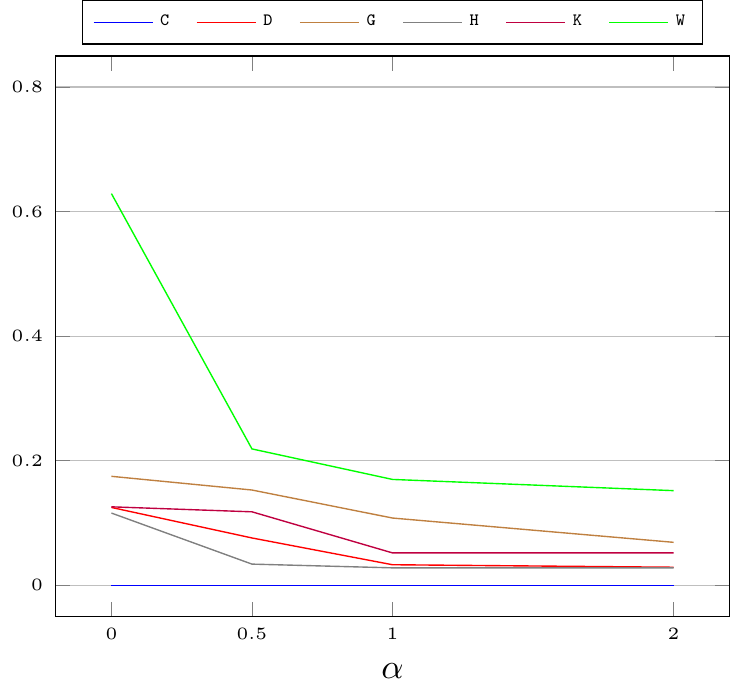}	
			\caption{$n=45$.}
		\end{subfigure}~
		\begin{subfigure}[b]{.45\linewidth}
			\includegraphics[scale=1]{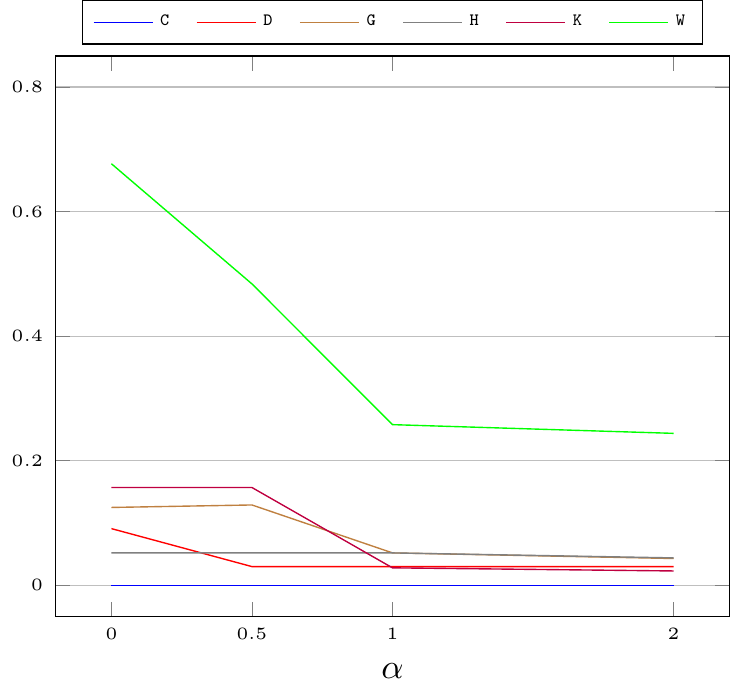}	
			\caption{$n=90$.}
		\end{subfigure}
		
		\begin{subfigure}[b]{.45\linewidth}
			\includegraphics[scale=1]{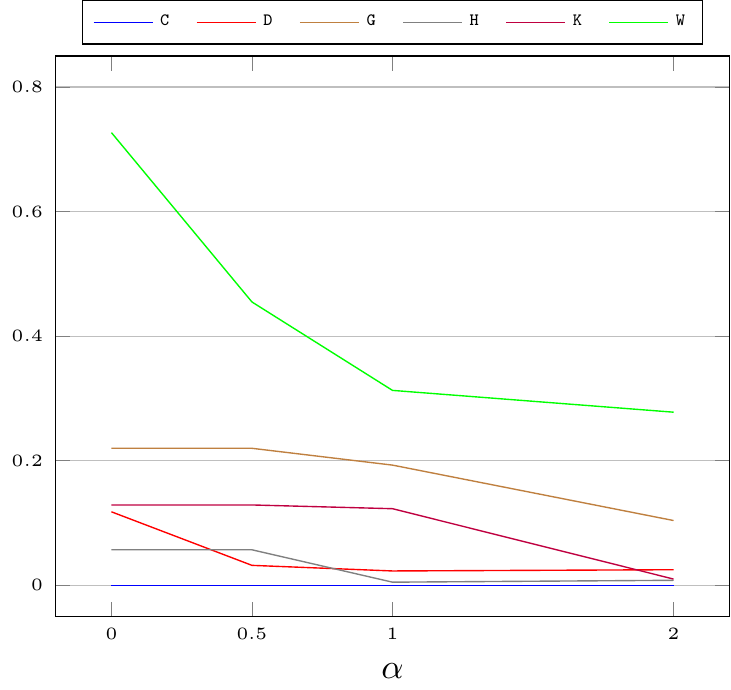}	
			\caption{$n=120$.}
		\end{subfigure}~
		\begin{subfigure}[b]{.45\linewidth}
			\includegraphics[scale=1]{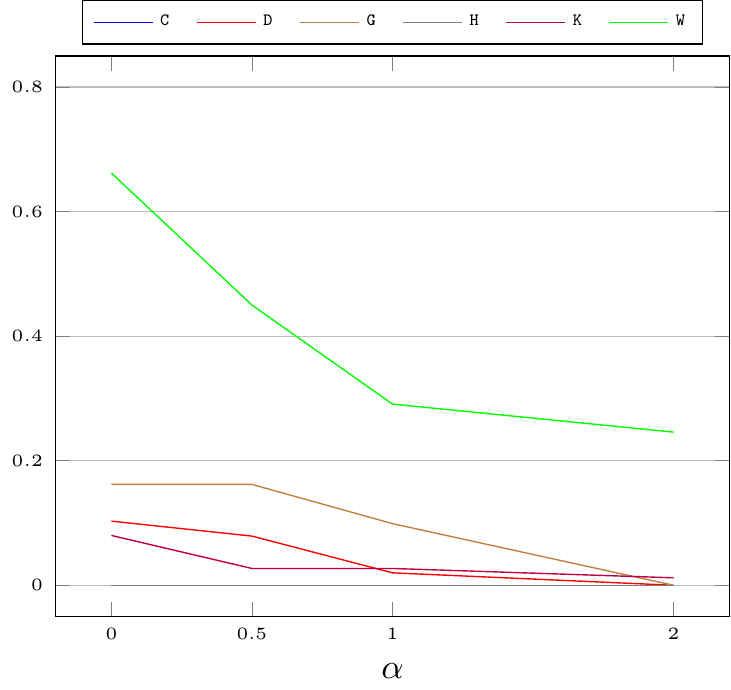}	
			\caption{$n=179$.}
		\end{subfigure}~
	\end{center}
	\caption{Price of efficiency averaged by $\alpha \in \{0,0.5,1,2\}$ and $n \in\{45,90,120,179\}$ for \texttt{Continuous} location problem.}
	\label{f:c-PoE-n}
\end{figure}

\begin{figure}[H]
	\begin{center}
		\begin{subfigure}[b]{.45\linewidth}
			\includegraphics[scale=1]{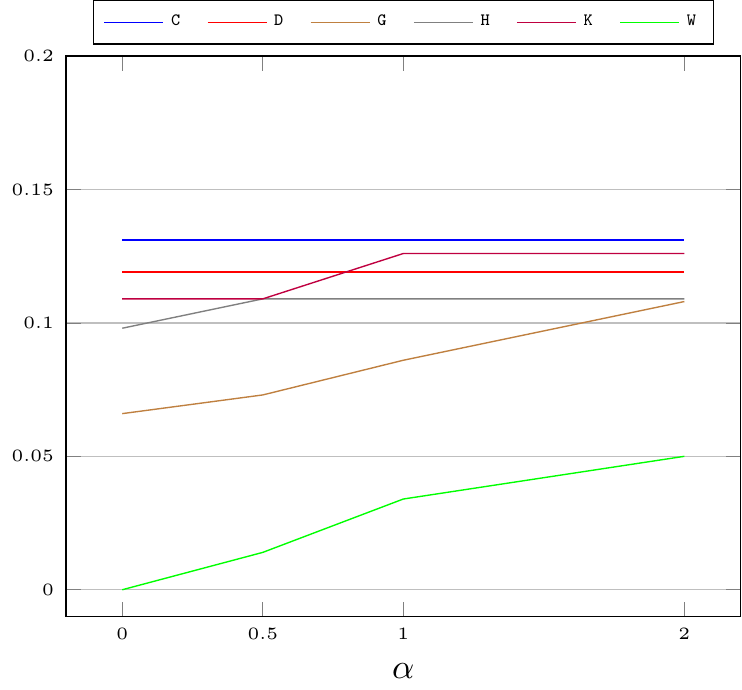}	
			\caption{$p=5$.}
		\end{subfigure}~
		\begin{subfigure}[b]{.45\linewidth}
			\includegraphics[scale=1]{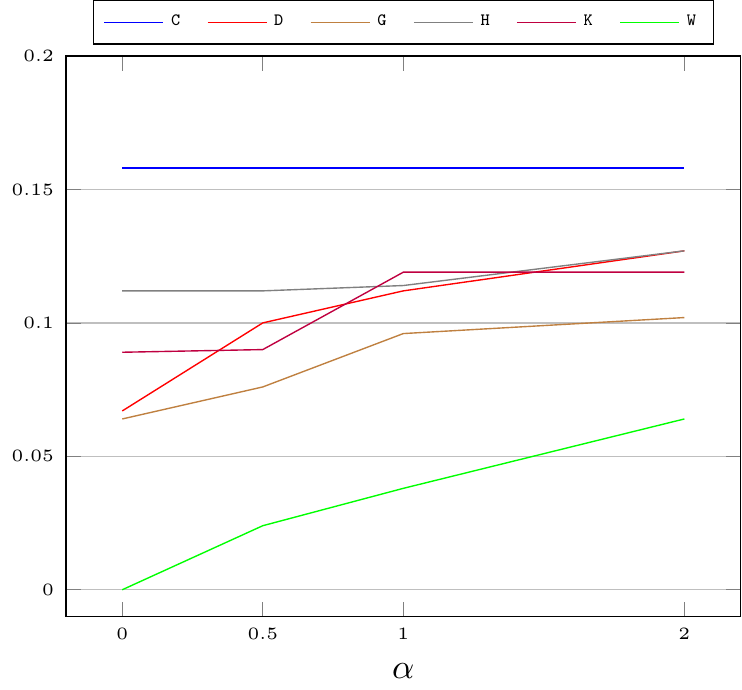}	
			\caption{$p=10$.}
		\end{subfigure}
		
		\begin{subfigure}[b]{.45\linewidth}
			\includegraphics[scale=1]{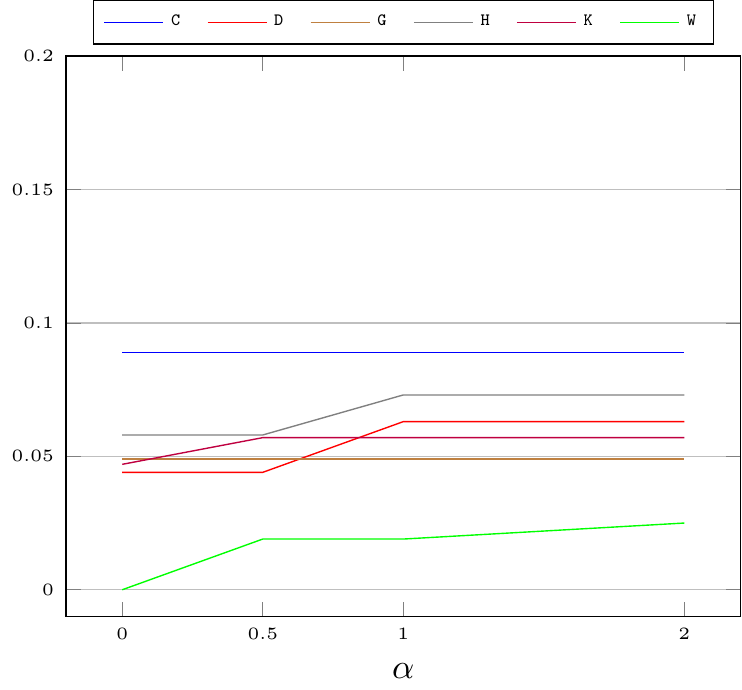}	
			\caption{$p=15$.}
		\end{subfigure}~
		\begin{subfigure}[b]{.45\linewidth}
			\includegraphics[scale=1]{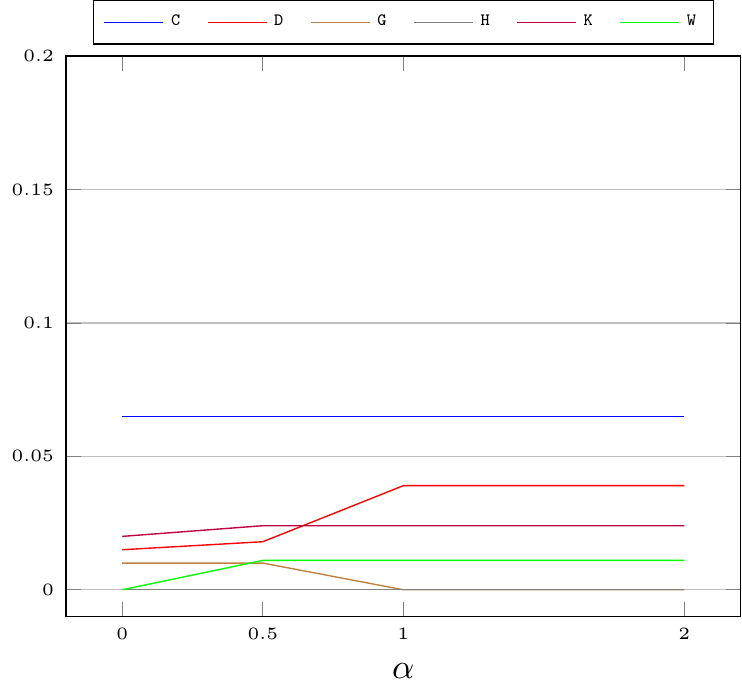}	
			\caption{$p=20$.}
		\end{subfigure}
	\end{center}
	\caption{Price of fairness averaged by $\alpha \in \{0,0.5,1,2\}$ and $p \in\{5,10,15,20\}$ \texttt{Continuous} location problem.}
	\label{f:c-PoF-p}
\end{figure}

\begin{figure}[H]
	\begin{center}
		\begin{subfigure}[b]{.45\linewidth}
			\includegraphics[scale=1]{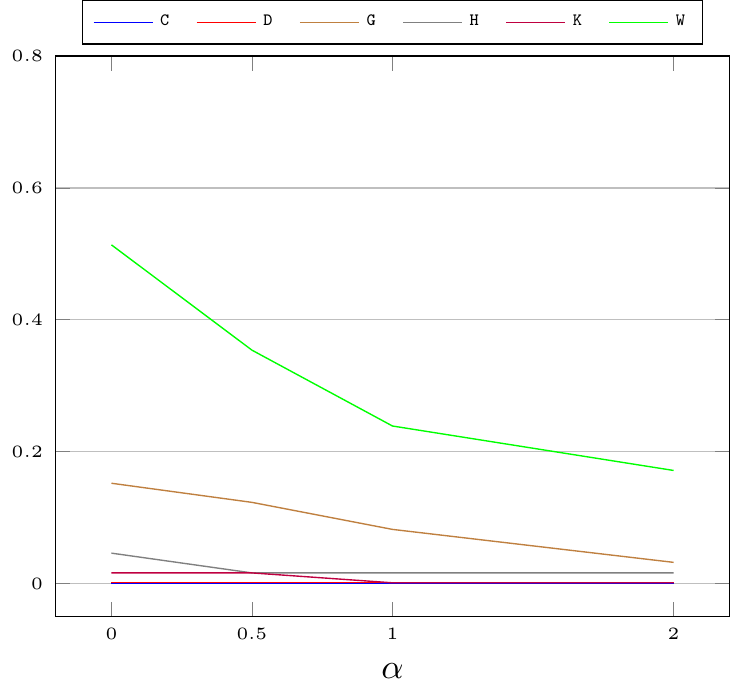}	
			\caption{$p=5$.}
		\end{subfigure}~
		\begin{subfigure}[b]{.45\linewidth}
			\includegraphics[scale=1]{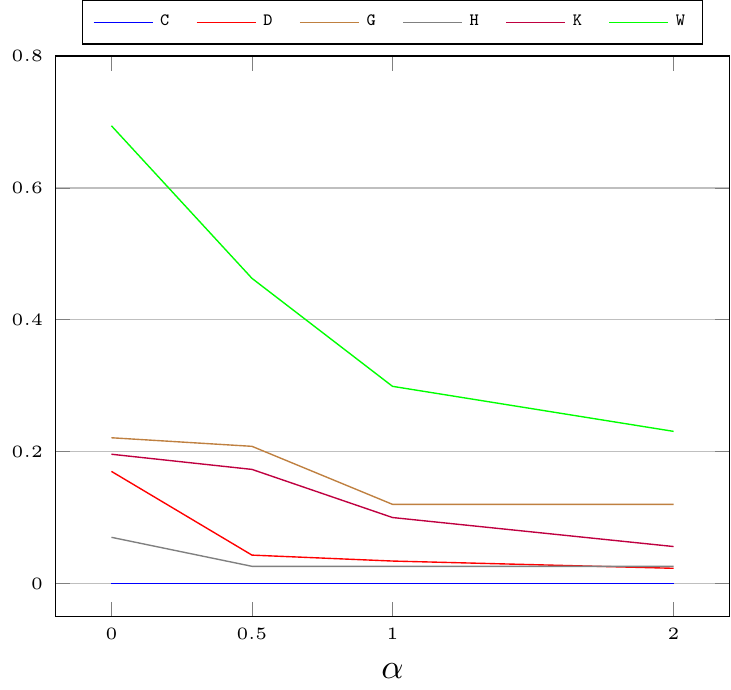}	
			\caption{$p=10$.}
		\end{subfigure}
		
		\begin{subfigure}[b]{.45\linewidth}
			\includegraphics[scale=1]{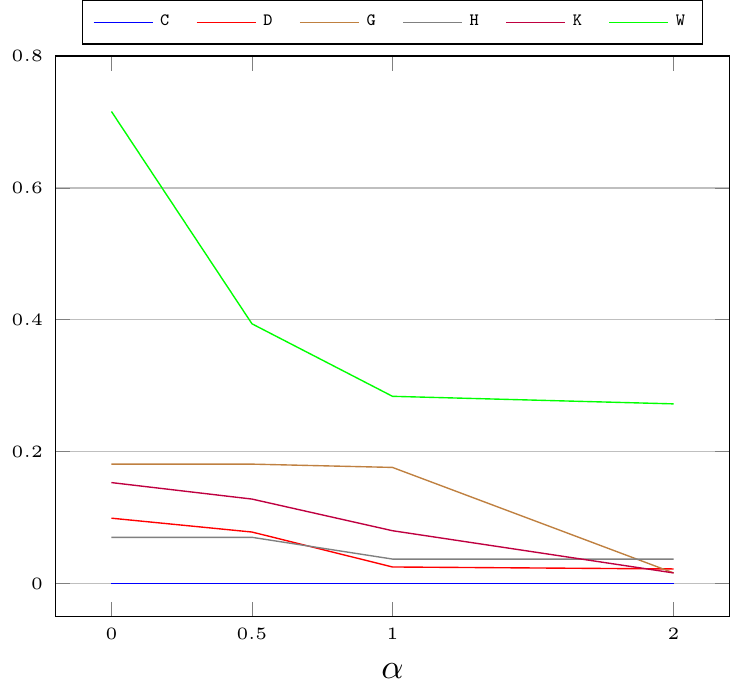}	
			\caption{$p=15$.}
		\end{subfigure}~
		\begin{subfigure}[b]{.45\linewidth}
			\includegraphics[scale=1]{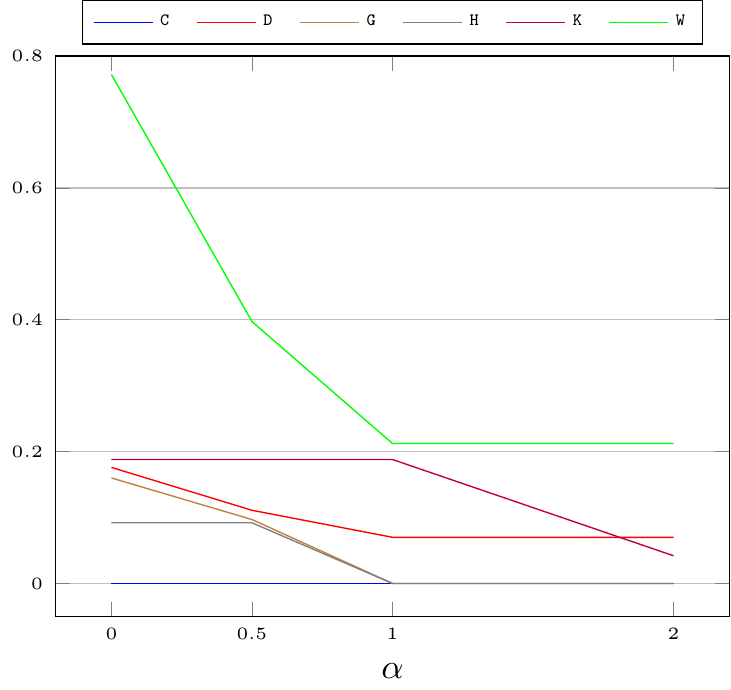}	
			\caption{$p=20$.}
		\end{subfigure}
	\end{center}
	\caption{Price of efficiency averaged by $\alpha \in \{0,0.5,1,2\}$ and $p \in\{5,10,15,20\}$ for \texttt{Continuous} location problem.}
	\label{f:c-PoE-p}
\end{figure}

\bibliographystyle{elsarticle-harv} 
\bibliography{faircovering}

\end{document}